\documentclass[a4paper]{article}

\usepackage{a4}
\usepackage{amsmath, amssymb}   %
\usepackage{mathtools}
\usepackage{graphicx,color}                   %
\usepackage{float}  
\usepackage{xspace}         
\usepackage[scientific-notation=true]{siunitx}         

\usepackage{tikz}
\tikzstyle{startstop} = [rectangle, rounded corners,minimum width=3cm, minimum height=1cm, text centered, text width=5cm, draw=black]
\tikzstyle{process} = [rectangle, minimum width=3cm, minimum height=1cm, text width=6 cm,text centered, draw=black]
\tikzstyle{decision} = [diamond, minimum width=3cm, minimum height=1cm, text width=3cm, text centered, draw=black]
\tikzstyle{arrow} = [thick,->,>=stealth]
\usetikzlibrary{spy}
\usetikzlibrary{shapes.geometric, arrows}
\usetikzlibrary{positioning} 
\usetikzlibrary{calc}
\usepackage{pgfplots}
\pgfplotsset{width=10cm,compat=1.9}
\usepackage[a4paper, total={6.5in, 8.5in}]{geometry}

\usepackage{hyperref}
\usepackage[all]{hypcap}

\definecolor{light_gray}{gray}{0.75}
\definecolor{lighter_gray}{gray}{0.5}
\colorlet{light_blue}{blue!20}
\definecolor{dark_green}{rgb}{0.0, 0.6, 0.0}
\definecolor{royal_blue}{rgb}{0.0, 0.22, 0.66}
\definecolor{salmon}{rgb}{1.0, 0.55, 0.41}
\definecolor{gold}{rgb}{0.8, 0.63, 0.21}
\definecolor{navy_blue}{rgb}{0.0, 0.0, 0.5}
\definecolor{crimson}{rgb}{0.79, 0.0, 0.09}
\definecolor{amethyst}{rgb}{0.6, 0.4, 0.8}
\definecolor{alizarin}{rgb}{0.82, 0.1, 0.26}
\definecolor{amaranth}{rgb}{0.9, 0.17, 0.31}
\definecolor{azure}{rgb}{0.0, 0.5, 1.0}
\definecolor{canaryyellow}{rgb}{0.82, 0.41, 0.12}
\definecolor{carrotorange}{rgb}{0.8, 0.33, 0.0}
\definecolor{cadmiumgreen}{rgb}{0.0, 0.42, 0.24}
\definecolor{copper}{rgb}{0.72, 0.45, 0.2}
\definecolor{aqua}{rgb}{0.5, 1.0, 0.83}
\definecolor{awesome}{rgb}{1.0, 0.13, 0.32}
\definecolor{candyapplered}{rgb}{1.0, 0.03, 0.0}
\definecolor{caribbeangreen}{rgb}{0.0, 0.8, 0.6}
\definecolor{aliceblue}{rgb}{0.94, 0.97, 1.0}
\definecolor{babyblue}{rgb}{0.54, 0.81, 0.94}

\newcommand{\leb}{\mathrm{leb}}
\newcommand{\lgl}{\mathrm{lgl}}
\newcommand{\bn}{\boldsymbol{n}}
\newcommand{\bx}{\boldsymbol{x}}
\newcommand{\by}{\boldsymbol{y}}

\newcommand{\bs}{\boldsymbol{s}}
\newcommand{\bo}{\boldsymbol{0}}
\newcommand{\mk}{\mathtt{k}}
\newcommand{\bA}{\mathrm{\textbf{A}}}

\newcommand{\bF}{\mathrm{\textbf{F}}}

\newcommand{\mr}{\mathrm{r}}
\newcommand{\me}{\mathrm{e}}
\newcommand{\mB}{\mathrm{B}}
\newcommand{\sas}{\mathrm{SAS}}
\newcommand{\abso}{\mathrm{abs}}
\newcommand{\sol}{\mathrm{sol}}

\newcommand{\ses}{\mathrm{SES}}
\newcommand{\sass}{\mathrm{SAS-S}}
\newcommand{\sess}{\mathrm{SES-S}}
\newcommand{\ith}{\mathrm{th}}

\newcommand{\ti}{\mathtt{i}}
\newcommand{\tk}{\mathtt{k}}
\newcommand{\te}{\mathtt{e}}

\newcommand{\ssphere}{\mathbb{S}^2}
\newcommand{\vry}{\varrho_{\scriptscriptstyle{Y}}}
\newtheorem{remark}{Remark}

\title{Domain Decomposition Method for Poisson--Boltzmann Equations based on Solvent Excluded Surface}
\author{Abhinav Jha\footnote{Institute of Applied Analysis and Numerical Simulation, University of Stuttgart,  Pfaffenwaldring 57, 70569,Stuttgart,Germany, \texttt{abhinav.jha@ians.uni-stuttgart.de}}, 
Benjamin Stamm\footnote{Institute of Applied Analysis and Numerical Simulation, University of Stuttgart,  Pfaffenwaldring 57, 70569,Stuttgart,Germany, \texttt{benjamin.stamm@ians.uni-stuttgart.de}}}

\date{}

\begin{document}
\maketitle
\begin{abstract}
In this paper, we develop a domain decomposition method for the nonlinear Poisson-Boltzmann equation based on a solvent-excluded surface widely used in computational chemistry. The model relies on a nonlinear equation defined in $\mathbb{R}^3$ with a space-dependent dielectric permittivity and an ion-exclusion function that accounts for steric effects. Potential theory arguments transform the nonlinear equation into two coupled equations defined in a bounded domain. Then, the Schwarz decomposition method is used to formulate local problems by decomposing the cavity into overlapping balls and only solving a set of coupled sub-equations in each ball. The main novelty of the proposed method is the introduction of a hybrid linear-nonlinear solver used to solve the equation. A series of numerical experiments are presented to test the method and show the importance of the nonlinear model.
\end{abstract}

\textbf{Keywords:} Implicit Solvation Model,  Poisson-Boltzmann Equation,  Domain Decomposition Method, Solvent Excluded Surface, Stern Layer

\section{Introduction}\label{sec:intro}
In computational chemistry, the Poisson-Boltzmann (PB) equation is a widely used model for modeling the ionic effects on molecular systems. It belongs to the class of implicit solvation models where the solute is treated microscopically, and the solvent is treated on its macroscopic physical properties, such as dielectric permittivity and ionic strength. Because of this treatment, implicit solvation models are computationally efficient, require fewer parameters, and implicitly consider the sampling over degrees of freedom of the solvent. For this reason, they are widely used in practice and are a popular computational approach to characterize solvent effects in the simulation of properties and processes of molecular systems \cite{TP94,RS99,OL00,TMC05}.

The history of the PB model can be traced back to 1910 when Gouy \cite{GS10} and Chapman \cite{CL13} independently used it to equate the chemical potential and relative forces acting on a small adjacent volume in an ionic solution between two plates having different voltages. Later, in 1923, Debye and H\"uckel generalized the concept by applying it to the theory of ionic solutions, leading to a successful interpretation of thermodynamic data  \cite{DH23}. A combination of both approaches, including a rigid layer close to the charged surface called the Stern layer and the Gouy-Chapman type diffusive layer, was introduced in 1924 by Stern \cite{Stern24}. The PB equation we consider in this paper is the realization of the Gouy-Chapman model with the possibility of including Stern layer correction \cite{ROHMR16}.

In this paper, we consider the nonlinear Poisson-Boltzmann (NPB) equation, which is used to describe the dimensionless electrostatic potential $\psi(\bx)$ and is given by
$$
-\nabla \cdot \left[\varepsilon(\bx)\nabla \psi(\bx) \right]+\lambda(\bx)\kappa^2\varepsilon_s\sinh\left(\psi(\bx)\right)=\frac{1}{\beta \varepsilon_{\abso}} \rho^{\sol}(\bx)\qquad \mathrm{in}\quad \mathbb{R}^3,
$$
for the case of $1:1$ electrolyte solvents. Here $\varepsilon(\bx)$ represents the relative space-dependent dielectric permittivity, $\varepsilon_\abso$ is the absolute dielectric permittivity of vacuum, $\varepsilon_s$ is the relative dielectric permittivity of the solvent, $\rho^{\mathrm{sol}}(\bx)$ is the charge distribution of the solute, $\lambda(\bx)$ is the ion-exclusion function that ensures that the ion concentration tends to zero inside the solute cavity, one in the bulk solvent region, and also accounts for the Stern layer, $\kappa$ is the Debye H\"uckel screening constant, and  $\beta = K_{\mathrm{B}}T/e$,  where $K_{\mathrm{B}}$ is the Boltzmann constant, $T$ is the temperature in Kelvins (K), and $e$ is the elementary charge.

With the assumption that $\psi$ satisfies the low-potential condition, i.e., $|\psi|\ll 1$, the NPB equation can be linearized to get the linear Poisson-Boltzmann (LPB) equation

$$
-\nabla \cdot \left[\varepsilon(\bx)\nabla \psi(\bx) \right]+\lambda(\bx)\kappa^2\varepsilon_s\psi(\bx)=\frac{1}{\beta \varepsilon_{\abso}} \rho^{\sol}(\bx)\qquad \mathrm{in}\quad \mathbb{R}^3.
$$

Some standard numerical methods for solving the PB equations include the finite difference method (FDM), the boundary element method (BEM), and the finite element method (FEM). We briefly overview them and mention some ongoing work in these areas. We refer to \cite{LZHM08} for a detailed comparison of the methods.

The FDM is one of the most popular methods for solving the LPB or the NPB equation. It follows the standard finite difference approach, where a grid covers the region of interest, and then different boundary conditions are chosen. Some of the popular software packages using the FDM include UHBD \cite{Madura95}, Delphi \cite{LiLi12}, MIBPB \cite{CCC10}, and APBS \cite{BSJ01, JE17}. One drawback of this approach is the increased computational cost with respect to the grid dimension, making it challenging to achieve high accuracy. Some recent developments in this area can be found in \cite{MTG10, MTH13, CEBG22}.

The BEM is another approach where the LPB equation is recast as an integral equation defined on a two-dimensional solute-solvent interface. This method can be optimized using the fast multipole or hierarchical treecode technique. The AFMPB solver \cite{LCH10, ZPH15} uses the former acceleration technique, whereas the TABI-PB \cite{WK13, WGK22} uses the latter one. 
The PB-SAM solver developed by Head-Gordon et al. \cite{LH06,YH10,YH13} discretizes the solute-solvent interface (such as the van der Waals (vdW) surface) with grid points on atomic spheres like a collocation method and solves the associated linear system by use of the fast multipole method. However, one drawback is that integral equation-based methods cannot be generalized to NPB. Hybrid approaches combining the FDM and BEM exist; see \cite{BF04}.

The FEM approach is one of the most flexible approaches for solving the PB equation. It can solve both the linear and the nonlinear PB, providing more flexible mesh refinement and proper convergence analysis \cite{CHX07}. A posteriori error estimation also exists for this method \cite{KNR20, NSRK21}. The SDPBS and SMPBS offer fast and efficient approximations of the size-modified PB equation \cite{Xie14, YX15, JXYXY15, XYX17}. Recently a hybrid approach combining the FEM and BEM has been proposed in \cite{BSBBC23}.

Now, we briefly overview the domain decomposition methods recently proposed in the context of implicit solvation models. Recently, in \cite{QSM19}, a domain decomposition algorithm for the LPB equation (ddLPB), which uses a particular Schwarz domain decomposition method, has been developed. The ddLPB method has been developed for the LPB equation defined on the vdW-cavity with discontinuous $\varepsilon(\bx)$ and didn't include a Stern layer correction, i.e., $\lambda(\bx)$. A further linear scaling approach for computing the first derivatives and eventually the forces has been presented in \cite{Jha23} following the ideas from \cite{MNS22}. The ideas of the ddLPB method can be traced back to the domain decomposition methods proposed for the conductor-like screening model (COSMO), (ddCOSMO) \cite{CMS13, LSCMM13, LSLS14} and the polarizable continuum model (PCM), (ddPCM) \cite{SCLM16, GLS17, NSSL19}. These methods do not require any mesh or grid of the molecular surface, are easy to implement, and ddCOSMO is about two orders of magnitude faster than the state-of-the-art \cite{LLSS14}.
An open-source software \texttt{ddX} has been released, which encompasses all these methods \cite{ddX}.
In general, domain decomposition methods for diffusion problems are known not to be scalable unless a coarse correction is used. However, for non-globular structures such as proteins, recent investigations \cite{CG18,CHS20,RS21} have shown the scalability of ddCOSMO without coarse correction.

An essential feature of the implicit solvation model is the choice of the solute-solvent boundary. Most methods use the vdW-cavity or the solvent-accessible surface (SAS), \cite{LR71} as they are topologically simple, but they do not describe the solute-solvent interaction well. The solvent-excluded surface (SES)  first developed in \cite{Ric77} is one of the few surfaces that captures the interaction quite well. In fact, it has also been observed that in some chemical calculations, the SES cavity yields better results \cite{RTL01, LMI17}. In \cite{QS16}, a mathematical framework was provided for computing the SES. As mentioned, the ddLPB method was developed with the vdW-cavity and can be extended to the SAS cavity. The ddPCM method with the SES cavity (ddPCM-SES) was proposed and studied in \cite{QSM18}.

This work proposes a new domain decomposition method for the NPB equation based on the SES (ddNPB-SES). We present a Schwarz domain decomposition method for the NPB equation that includes steric effects, i.e., the presence of a Stern layer, which was missing from the previous work, \cite{QSM19}. A further development from the previously mentioned work is using the SES for the solute-solvent boundary. With the help of the SES, we introduce a continuous relative dielectric permittivity function $\varepsilon(\bx)$ and an ion-exclusion function $\lambda(\bx)$. Based on the definition of $\varepsilon(\bx)$ and $\lambda(\bx)$, we develop a hybrid linear-nonlinear solver that uses local spectral methods for the discretization, using spherical harmonics and Legendre polynomials as basis functions. This contrasts with the traditional numerical techniques for solving the NPB equation, where a nonlinear solver is used in the whole domain. The development of a hybrid linear-nonlinear solver becomes important as the low-potential condition is usually not satisfied near the surface of the molecules. Hence, we use a nonlinear model near the surface of a molecule and a linear model in the bulk solvent region. The nonlinearity of the model also plays a critical role in the computation of the electrostatic solvation energy $E_s$, as under the low-potential condition, the solvent's contribution is not considered. The proposed method can also solve the LPB equation defined on the SES, i.e., ddLPB-SES; hence, this work encompasses both models. Therefore, this work is an extension of the ddLPB method, \cite{QSM19}, and introduces a new method, ddNPB-SES, for treating the important nonlinear effects. The major contribution of this work can be summarized as:
\begin{enumerate}
\item[$\bullet$] Introduction of a domain decomposition method for the NPB equation.
\item[$\bullet$] Use of a Solvent-Excluded Surface for the solute-solvent boundary.
\item[$\bullet$]  Inclusion of steric effects.
\end{enumerate}
To the best of our knowledge, this is the first work that deals with all three points mentioned above. This work can also be regarded as a PDE-IE (Integral Equation) coupling of the NPB equation. 
In this regard, the recent work \cite{BSBBC23} is related as it deals with such a coupling, but it considers the LPB equation only.

Finally, let us mention that the method's efficiency is based on the scalability of Poisson and screened Poisson problems in the cavity of the molecule. 
Therefore, the scalability discussed above also applies in this setting.

The breakdown of the paper is as follows: In Sec.~\ref{sec:prob_statement}, we derive the NPB equation and introduce different solute-solvent boundaries. We also present a continuous relative dielectric permittivity function and an ion-exclusion function based on the SES. In Sec.~\ref{sec:form_prob}, we transform the problem into different domains, introduce a global strategy for solving them, and lay out the domain decomposition method. In Sec.~\ref{sec:single_domain}, we derive single-domain solvers for the homogeneous screened Poisson (HSP) and the generalized screened Poisson (GSP) equation in a unit ball. Next, in Sec.~\ref{sec:numres}, we present a comprehensive numerical study for the ddNPB-SES method for molecules ranging from one to twenty-four atoms. Lastly, in Sec.~\ref{sec:summary}, we present a summary and an outlook.

\section{Problem Statement}\label{sec:prob_statement}
We represent the solvent by a polarizable and ionic continuum. The freedom of the movement of ions is modeled by Boltzmann statistics, i.e., the Boltzmann equation is used to calculate the local ion density $c_i$ of the $i^{\ith}$ type of ion as follows
\begin{equation}\label{eq:ion_density}
c_i=c_i^{\infty}\exp\left( \frac{-W_i}{K_{\mathrm{B}}T}\right),
\end{equation}
where $c_i^{\infty}$ is the bulk ion concentration at an infinite distance from the solute molecule, and $W_i$ is the work required to move the $i^\ith$ type of ion to a given position from an infinite distance.

The electrostatic potential $\tilde{\psi}(\bx)$ of a general implicit solvation model is described by the Poisson equation as follows

\begin{equation}\label{eq:poisson_eqn}
-\nabla\cdot \left[\varepsilon_\abso\varepsilon(\bx)\nabla \tilde{\psi}(\bx)\right] =  \rho^{\mathrm{sol}}(\bx) + \rho^{\mathrm{ions}}(\bx)\quad \mathrm{in}\ \ \mathbb{R}^3,
\end{equation}
where $\tilde{\psi}(\bx) = \mathcal{O}\left(1/|\bx|\right)$ as $|\bx|\rightarrow \infty$. Here $\rho^{\mathrm{ions}}(\bx)$ is the charge distribution of the solvation system.

We derive the PB equation using Eq.~\eqref{eq:ion_density} and Eq.~\eqref{eq:poisson_eqn} as
\begin{equation}\label{eq:pb_eqn}
-\nabla \cdot \left[\varepsilon_\abso\varepsilon(\bx)\nabla \tilde{\psi}(\bx)\right] = \rho^{\mathrm{sol}}(\bx) + \lambda(\bx)\sum_{i=1}^{N_{\mathrm{ions}}}c_i^{\infty}z_ie\exp\left( \frac{-z_ie\tilde{\psi}(\bx)}{K_{\mathrm{B}}T}\right)\quad \mathrm{in}\ \ \mathbb{R}^3,
\end{equation}
where $z_i$ is the (partial) charge of the $i^{\ith}$ type of ion and $N_{\mathrm{ions}}$ are the number of ion types.

In the case of a $1:1$ electrolyte, there are two opposite charge ions $(+e$ and $-e)$ and we then get
\begin{alignat}{3}\label{eq:derive_nlpb}
\sum_{i=1}^2c_i^{\infty}z_ie\exp\left( \frac{-z_ie\tilde{\psi}(\bx)}{K_{\mathrm{B}}T}\right) & =  ce\exp\left(\frac{-e\tilde{\psi}(\bx)}{K_{\mathrm{B}}T}\right) - ce\exp\left(\frac{e\tilde{\psi}(\bx)}{K_{\mathrm{B}}T}\right)\nonumber \\
& =  -2ce\sinh\left(\frac{e\tilde{\psi}(\bx)}{K_{\mathrm{B}}T}\right).
\end{alignat}

Then substituting Eq.~\eqref{eq:derive_nlpb} into Eq.~\eqref{eq:pb_eqn} we obtain
\begin{equation}\label{eq:npb_eqn}
-\nabla \cdot \left[\varepsilon_\abso\varepsilon(\bx)\nabla \tilde{\psi}(\bx)\right] +\lambda(\bx)2 ce\sinh\left(\frac{e\tilde{\psi}(\bx)}{K_{\mathrm{B}}T}\right)= \rho^{\mathrm{sol}}(\bx) \quad \mathrm{in}\ \ \mathbb{R}^3,
\end{equation}
the nonlinear Poisson-Boltzmann equation (NPB). 

\subsection{Solute Probe}\label{sec:solute_probe}
One of the important properties of the implicit solvation model is the choice of the solute probe and, accordingly, the solute-solvent boundary. A straightforward choice is using the van der Waals (vdW) surface, i.e., the topological boundary of the union of solute's vdW-atoms with experimentally fitted radii. Another choice is the solvent-accessible surface (SAS), which is defined by tracing the center of an idealized (spherical) solvent probe (representing a solvent molecule) when rolling over the solute molecule. The region enclosed by the SAS is called the SAS-cavity, which we denote by $\Omega_\sas$ and its boundary by $\Gamma_\sas$.

The vdW and the SAS are topologically not the correct answers to the cavity problem as they poorly describe the region where the solvent can touch. However, as they are topologically simple, they are widely used in numerical computations. Another solute-solvent boundary is the solvent-excluded surface (SES), which represents the boundary of the region where the probe has no access due to the presence of the solute. In other words, SES is the boundary of the union of all spherical probes that do not intersect the vdW balls of the solute molecule. The region enclosed by the SES is the SES cavity, which we denote by $\Omega_{\ses}$ and the boundary by $\Gamma_\ses$. The mathematical characterization of the surface can be found in \cite{QS16}. From the geometrical viewpoint, the SES is smoother than the vdW and SAS cavities.

The NPB equation does not consider the absorbing ions' finite size; hence, the ionic concentration can exceed the maximally allowed coverage near the surface. These are referred to as steric effects. The NPB model is modified to include a Stern layer, \cite{BAO97} to account for these steric effects.

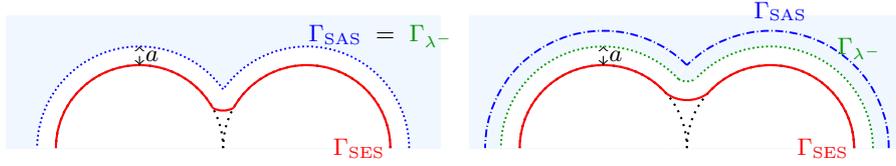
\begin{figure}[t!]
\centering
\begin{tikzpicture}[scale=1.1]
\fill[line width=0.8pt, color = aliceblue] (-1.6,0) -- (3.6,0) -- (3.6,1.6) -- (-1.6,1.6) -- cycle;
\filldraw[white] (1.22474487139,0) arc(0: 180:1.22474487139);
\filldraw[white] (3.22474487139,0) arc(0: 180:1.22474487139);
\draw[line width=0.8pt,black, dotted] (1,0) arc(0: 180:1);
\draw[line width=0.8pt,red] (0.875, 0.4841) arc(28.9538754: 180:1);
\draw[line width=0.8pt,black, dotted] (3,0) arc(0: 180:1);
\draw[line width=0.8pt,red] (3,0) arc(0: 151.0461246:1);
\filldraw[white] (0.875,.4841) arc(240.7277123: 299.2722877:0.25564428411);
\draw[line width=0.8pt,red] (0.875,.4841) arc(240.7277123: 299.2722877:0.25564428411);
\draw[densely dotted, line width=0.7pt,blue] (1, 0.7071) arc(35.2641307: 180:1.22474487139);
\draw[densely dotted, line width=0.7pt,blue] (3.225,0) arc(0:144.7358693:1.22474487139);
\draw [dashed, line width=0.4pt, to-to](0,1) -- (0,1.2247);
\begin{scriptsize}
\draw (3,0) node[ text width= 1.5cm, align=left] {\small{$\textcolor{red}{\Gamma_{\ses}}$}};
\draw (2.62,1.35) node[ text width= 2.5cm, align=right] {\small{$\textcolor{blue}{\Gamma_{\sas}} = \textcolor{dark_green}{\Gamma_{\lambda^-}}$}};
\draw (0.15,1.1135) node[align=right] {\small{$a$}};
\end{scriptsize}
\end{tikzpicture}\hspace*{0.001em}
\begin{tikzpicture}[scale=1.1]
\fill[line width=0.8pt, color = aliceblue] (-1.6,0) -- (3.6,0) -- (3.6,1.6) -- (-1.6,1.6) -- cycle;
\filldraw[white] (1.22474487139,0) arc(0: 180:1.22474487139);
\filldraw[white] (3.22474487139,0) arc(0: 180:1.22474487139);
\filldraw[white] (0.7511, 0.6602) arc(233.5925939: 306.4074061:0.42121372247)--(1,1)--cycle;
\filldraw[aliceblue] (0.8995, 0.8312) arc(239.2317219:300.7682781:0.19645512861)--(1,1)--cycle;
\draw[line width=0.8pt,black, dotted] (1,0) arc(0: 180:1);
\draw[line width=0.8pt,red] (0.75, 0.6614) arc(41.4079965: 180:1);
\draw[line width=0.8pt,black, dotted] (3,0) arc(0: 180:1);
\draw[line width=0.8pt,red] (3,0) arc(0: 138.5920035:1);
\draw[line width=0.8pt,red] (0.7511, 0.6602) arc(233.5925939: 306.4074061:0.42121372247);
\draw[densely dotted, line width=0.7pt,dark_green] (0.8995, 0.8312) arc(42.7400635: 180:1.22474487139);
\draw[densely dotted, line width=0.7pt,dark_green] (3.225,0) arc(0:137.2599365:1.22474487139);
\draw[densely dotted, line width=0.7pt,dark_green] (0.8995, 0.8312) arc(239.2317219:300.7682781:0.19645512861);
\draw[densely dashdotted, line width=0.7pt,blue] (1,1) arc(45: 180:1.414);
\draw[densely dashdotted, line width=0.7pt,blue] (3.414,0) arc(0: 135:1.414);
\draw [dashed, line width=0.4pt, to-to](0,1) -- (0,1.2247);
\begin{scriptsize}
\draw (3,0) node[ text width= 1.5cm, align=left] {\small{$\textcolor{red}{\Gamma_{\ses}}$}};
\draw (2.65,1.2) node[ text width= 1.5cm, align=right] {\small{$\textcolor{dark_green}{\Gamma_{\lambda^-}}$}};
\draw (1.75,1.65) node[ text width= 1.5cm, align=right] {\small{$\textcolor{blue}{\Gamma_{\sas}}$}};
\draw (0.15,1.1135) node[align=right] {\small{$a$}};
\end{scriptsize}
\end{tikzpicture}
\caption{SES surface with a probe of radius $r_p$ (left) and with $r_p+a$ (right) .}\label{fig:ses_and_ses_a}
\end{figure}

Now, we set certain notations. We assume that the molecule comprises $M$ atoms and the $i^\ith$ atom has center $\bx_i$ and vdW radii $r_i$. The solvent probe radius is denoted by $r_p$, and the Stern layer length is denoted by $a$. In practice, the SAS cavity is defined by rolling a probe of radius $r_p$; to incorporate the Stern layer instead, we use a probe of radius $r_p+a$. Including the Stern layer means the solute touches the solvent after distance $a$. If we use a probe of radius $r_p$, we have a non-smooth surface near the solute-solvent interaction (see Fig.~\ref{fig:ses_and_ses_a} (left)). Therefore, we use a probe of radius $r_p+a$ to restore smoothness near the interaction (see Fig.~\ref{fig:ses_and_ses_a} (right)). The boundary $\Gamma_{\lambda^-}$ denotes the boundary where the ion concentration starts. More details will be given in the following subsection. Furthermore, for each atom, we define an ``enlarged" ball $\Omega_i$ with center $\bx_i$ and radius $R_i=r_i+r_p+a+r_0$, where $r_0$ is a non-negative constant used to control the nonlinear regime. Outside of the enlarged balls, one can use the LPB equation as the low potential condition is satisfied. The parameter $r_0$ gives a good idea for the cutoff where we switch from the nonlinear model to the linear model.

The SES cavity is entirely covered by the union of $\Omega_0$ of enlarged balls, i.e.,
$$
\Omega_\ses\subset \Omega_0:=\bigcup_{i=1}^M\Omega_i,\quad \mathrm{where}\quad \Omega_i=B_{R_i}(\bx_i),
$$
and $B_{R_i}(\bx_i)$ denotes a ball with center $\bx_i$ and radius $R_i$. We also denote the bulk solvent region by $\Omega_{\infty}:=\mathbb{R}^3\setminus\Omega_0=\Omega_0^{\mathrm{C}}$ and the boundary of $\Omega_0$ by $\Gamma_0$.

Let $f_\sas$ denote the signed distance function to $\Omega_{\sas}$ (i.e., negative inside the SAS cavity and positive outside the SAS cavity). We then have a mathematical characterization of the two cavities
$$
\Omega_{\ses} = \left\lbrace \bx\in \mathbb{R}^3:\ f_\sas(\bx)\leq -r_p-a\right\rbrace\ \  \mathrm{and}\ \  \Omega_0=\left\lbrace \bx\in \mathbb{R}^3:\ f_\sas(\bx)\leq r_0\right\rbrace.
$$
Also, we can characterize their boundary surfaces by
$$
\Gamma_\ses = f^{-1}_\sas(-r_p-a)\quad \mathrm{and}\quad \Gamma_0=f^{-1}_\sas(r_0).
$$

\begin{remark}\label{rem:rho_sol_def}
It is reasonable to assume that the solute's charge distribution $\rho^{\sol}$ is supported in $\Omega_0$. In this paper we consider the classic description of $\rho^{\sol}$ given by
$$
\rho^{\sol}(\bx)=\sum_{i=1}^Mq_i\delta(\bx-\bx_i),
$$
where $q_i$ denotes the (partial)\ charge carried by the $i^{\ith}$ atom, and $\delta$ is the Dirac-delta function. For a quantum description of the solute, $\rho^{\sol}$ comprises a sum of nuclear charges and the electron charge density.
\end{remark}

\subsection{Dielectric Permittivity and Ion-Exclusion Function}\label{sec:ion_exclusion}
In this subsection, we construct an SES-based dielectric permittivity and an ion-exclusion function associated with $f_\sas$, which follows the ideas from \cite{QSM18}.

The solvent dielectric permittivity is assumed to be constant and equal to the bulk dielectric permittivity after a certain distance from the solute molecule. This is a reasonable assumption as the solvent density at positions far from the solute molecule reaches the bulk limit.

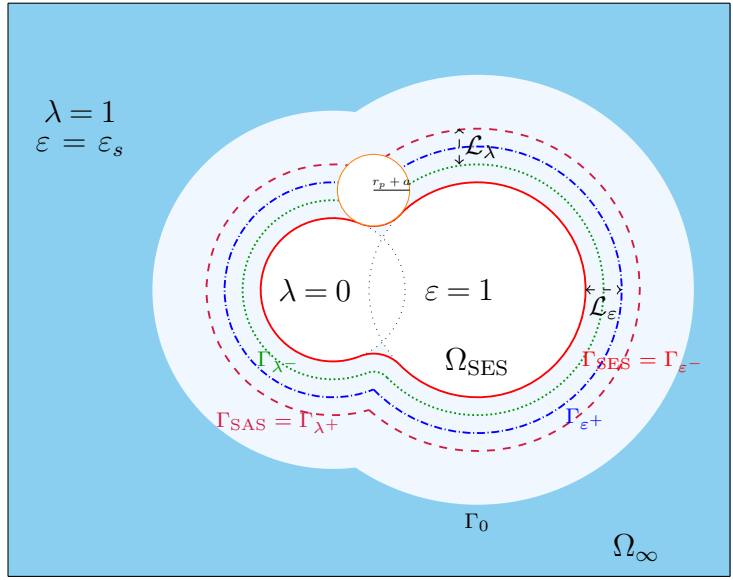
\begin{figure}[t!]
\centering
\begin{tikzpicture}[scale=0.95]
\fill[line width=0.8pt, color = babyblue] (1,0) -- (11,0) -- (11,8) -- (1,8) -- cycle;
\draw (1,0) -- (11,0);
\draw (11,0) -- (11,8);
\draw (11,8) -- (1,8);
\draw (1,8) -- (1,0);
\fill[line width=0.8pt, color = aliceblue] (5.5,4) circle [radius= 2.5];
\fill[line width=0.8pt, color = aliceblue] (7.5,4) circle [radius=3];
\fill[line width=0.8pt, color = white] (5.5,4) circle [radius=1];

\fill[line width=0.8pt, color = white	] (7.5,4) circle [radius=1.5];
\draw[dotted, black, line width = 0.3] (5.5,4) circle(1 cm);
\draw[dotted, black, line width = 0.3] (7.5,4) circle(1.5cm);

\draw[line width=0.7pt,red] (5.875338956423, 4.926887611996) arc(247.9756872: 321.9513743:0.505);
\draw[line width=0.7pt,red] (6.42187223012, 2.957099951182) arc(44.0486257: 112.0243128:0.5);
\draw[line width=0.7pt,red] (5.875338956423, 4.926887611996) arc(67.9547342: 292.0452658:1);
\draw[line width=0.7pt,red] (6.42187223012, 2.957099951182) arc(224.0484735:494.0614503:1.5);
\draw[densely dotted, line width=0.7pt,dark_green] (6.125	,5.082531754731) arc(60:290:1.25);
\draw[densely dotted, line width=0.7pt,dark_green] (6.25086911	,2.774368) arc(224.4559857:501.7867893:1.75);

\draw[densely dotted, line width=0.7pt,dark_green] (6.25086911	,2.774368) arc(41.871835:122.0447562:0.25);

\draw[densely dashdotted, line width=0.7pt,blue] (6.0625	,5.39053721633) arc(67.9756872:292.0243128:1.5);
\draw[densely dashdotted, line width=0.7pt,blue] (6.0625,2.60946278367) arc(224.0486257:495.9513743:2);

\draw[dashed, line width=0.7pt,alizarin] (6,	5.677050983125) arc(73.3984504 :286.6015496 :1.75);
\draw[dashed, line width=0.7pt,alizarin] (6,	2.322949016875) arc(228.1896851:491.8103149:2.25);
\filldraw[line width=0.3pt, color = orange, fill = white] (6.0625,5.390537) circle [radius=0.5];
\draw [dashed, line width=0.4pt, to-to](9,4) -- (9.5,4);
\draw [dashed, line width=0.4pt, to-to](7.25,5.75) -- (7.25,6.25);
\draw [line width=0.2pt](6.0625,5.390537) -- (6.5625,5.390537);
\begin{scriptsize}
\draw (5.25,4) node[ text width= 1.5cm, align=center] {\large{$\lambda$ = 0}};
\draw (2,6.5) node[ text width= 1.5cm, align=center] {\large{$\lambda$ = 1}};

\draw (7.25,4) node[ text width= 1.5cm, align=center] {\large{$\varepsilon$ = 1}};
\draw (2,6) node[ text width= 1.5cm, align=center] {\textcolor{black}{\large{$ \varepsilon= \varepsilon_s$}}};
\draw (7.5,3) node[ text width= 1.5cm, align=center] {\large{$\Omega_\ses$}};
\draw (9.7,0.4) node[ text width= 1.5cm, align=center] {\large{$\Omega_{\infty}$}};

\draw (9.8,3) node[ text width= 1.75cm, align=center] {\textcolor{red}{\footnotesize{$\Gamma_{\ses} =\Gamma_{\varepsilon^-}$}}};
\draw (4.75,3) node[ text width= 1.5cm, align=center] {\textcolor{dark_green}{\footnotesize{$\Gamma_{\lambda^-}$}}};
\draw (9,2.2) node[ text width= 1.5cm, align=center] {\textcolor{blue}{\footnotesize{$\Gamma_{\varepsilon^+}$}}};
\draw (4.75,2.15) node[ text width= 1.75cm, align=center] {\textcolor{alizarin}{\footnotesize{$\Gamma_\sas =\Gamma_{\lambda^+}$}}};
\draw (7.5,0.75) node[ text width= 1.5cm, align=center] {\footnotesize{$\Gamma_0$}};

\draw (9.25,3.75) node[ text width= 1.5cm, align=center] {\normalsize{$\mathcal{L}_{\varepsilon}$}};
\draw (6.3125,5.490537) node[ text width= 1cm, align=center] {\scalebox{.65}{$r_p+a$}};
\draw (7.55,6) node[ text width= 1.5cm, align=center] {\normalsize{$\mathcal{L}_{\lambda}$}};
\end{scriptsize}
\end{tikzpicture}
\caption{Solute probes and solute-solvent boundary for a diatmoic molecule.}\label{fig:solute_cavity}
\end{figure}

Taking the SES as the solute-solvent boundary implies that the dielectric permittivity in the SES cavity is always one, i.e.,  the relative dielectric permittivity of vacuum and constant outside $\sas$,  with value $\varepsilon_s$. Let us denote the layer where $\varepsilon(\bx)$ varies by $\mathcal{L}_{\varepsilon}:=\Omega_{\varepsilon^+}\setminus \Omega_{\varepsilon^-}$ where
$$
\Omega_{\varepsilon^-}=\Omega_{\ses}\ \  \mathrm{and}\ \ \Omega_{\varepsilon^+} = \left\lbrace \bx\in \mathbb{R}^3:\ f_\sas(\bx)\leq -a\right\rbrace,
$$
with $\Gamma_{\varepsilon^-}$ and $\Gamma_{\varepsilon^+}$ as their respective boundary. Similarly, the ion-exclusion function is zero inside $\Omega_{\lambda^-}$ and one outside $\Omega_{\lambda^+}$, where
$$
\Omega_{\lambda^-} = \left\lbrace \bx\in \mathbb{R}^3:\ f_\sas(\bx)\leq -r_p\right\rbrace\ \  \mathrm{and}\ \ \Omega_{\lambda^+} = \Omega_{\sas},
$$
with $\Gamma_{\lambda^-}$ and $\Gamma_{\lambda^+}$ as their respective boundary. We denote the layer where $\lambda(\bx)$ varies by $\mathcal{L}_\lambda$, i.e., $\mathcal{L}_\lambda = \Omega_{\lambda^+}\setminus \Omega_{\lambda^-}$. Fig.~\ref{fig:solute_cavity} represents all the molecular probes and the boundaries introduced.

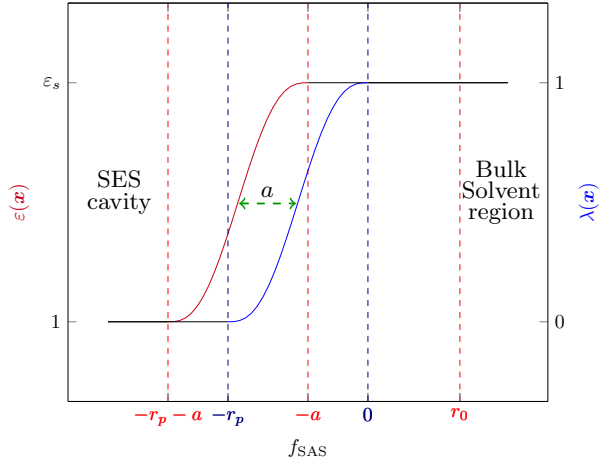
\begin{figure}
\centering
\begin{tikzpicture}[scale=0.75]
\begin{axis}[
  axis y line*=left,
  ymin=0, ymax=5,
  ytick = {1,4},
  yticklabels = {1, $\varepsilon_s$},
  xtick = {-2.5, -1, 1, 2.5, 4.8},
  xticklabels = {$\textcolor{red}{-r_p-a}$, $\textcolor{navy_blue}{-r_p}$, $\textcolor{red}{-a}$, $\textcolor{navy_blue}{0}$, $\textcolor{red}{r_0}$},
  ylabel=$\textcolor{crimson}{\varepsilon(\bx)}$,
]
\pgfplotsinvokeforeach{-2.5, 1,4.8}{
  \draw[dashed, color = red] ({rel axis cs: 0,0} -| {axis cs: #1, 0}) -- ({rel axis cs: 0,1} -| {axis cs: #1, 0});}
\addplot[domain=-4:-2.5] {1};
\addplot[domain=-2.5:1, color=crimson] {0.03427*x*x*x*x*x + 0.12852*x*x*x*x - 0.15708*x*x*x - 0.64259*x*x+1.07098*x + 3.5659};
\addplot[domain=1:6] {4};
\end{axis}
\begin{axis}[
  axis y line*=right,
  ymin=0, ymax=5,
  ytick = {1,4},
  yticklabels = {0, 1},
   xtick = {-2.5, -1, 1, 2.5, 4.8},
  xticklabels = {$\textcolor{red}{-r_p-a}$, $\textcolor{navy_blue}{-r_p}$, $\textcolor{red}{-a}$, $\textcolor{navy_blue}{0}$, $\textcolor{red}{r_0}$},
  xlabel=$f_{\sas}$,
  ylabel=$\textcolor{blue}{\lambda(\bx)}$,
]
\pgfplotsinvokeforeach{-1, 2.5}{
  \draw[dashed, color = navy_blue] ({rel axis cs: 0,0} -| {axis cs: #1, 0}) -- ({rel axis cs: 0,1} -| {axis cs: #1, 0});}
\addplot[domain=-4:-1] {1};
\addplot[domain=-1:2.5, color=blue] {  0.03427*x*x*x*x*x  -0.12852*x*x*x*x -0.15708*x*x*x + 0.64259*x*x+1.07098*x + 1.43410};
\addplot[domain=2.5:6] {4};
\end{axis}
\draw [dashed, line width=0.7pt, to-to, color = dark_green](3,3.5) -- (4,3.5);
\begin{scriptsize}
\draw (3.5,3.7) node[ text width= 1.5cm, align=center] {\small{$a$}};
\draw (0.9,3.7) node[ text width= 1.5cm, align=center] {\small{SES cavity}};
\draw (7.6,3.7) node[ text width= 1.5cm, align=center] {\small{Bulk Solvent region}};
\end{scriptsize}
\end{tikzpicture}
\caption{Schematic diagram of the dielectric permittivity, $\varepsilon(\bx)$ (left $y$-axis) and the ion-exclusion function, $\lambda(\bx)$ (right $y$-axis) with respect to $f_{\sas}$.}\label{fig:ion_exclusion}
\end{figure}

The remaining work is to determine $\varepsilon(\bx)$ and $\lambda(\bx)$ in the intermediate layer, $\mathcal{L}_\varepsilon$ and $\mathcal{L}_\lambda$, respectively.  We choose the following modified definition of the permittivity function from \cite{QSM18},
\begin{equation}
\varepsilon(\bx)=\left\lbrace 
\begin{array}{ll}
1 & \bx \in \Omega_\ses,\\
1+(\varepsilon_s-1)\xi\left(\dfrac{f_\sas(\bx)+r_p + a}{r_p}\right) & \bx \in \mathcal{L}_{\varepsilon},\\
\varepsilon_s & \mathrm{else},
\end{array}\right.
\end{equation}
and define the ion-exclusion function as
\begin{equation}
\lambda(\bx)=\left\lbrace 
\begin{array}{ll}
0 & \bx \in \Omega_{\lambda^-} ,\\
\xi\left(\dfrac{f_\sas(\bx)+ r_p}{r_p}\right) & \bx \in \mathcal{L}_{\lambda},\\
1 & \mathrm{else},
\end{array}\right.
\end{equation}
where $\xi(\cdot)$ is a continuous function defined on $[0,1]$, satisfying $\xi(0)=0$,\ $\xi(1)=1,\ \xi'(0)=0$, and $\xi'(1)=0$. $\varepsilon(\bx)$ and $\lambda(\bx)$ can be seen as distance-dependent functions where the ``distance" represents the signed distance to SAS, see Fig.~\ref{fig:ion_exclusion} for a schematic diagram. The function $\xi(\cdot)$ can be chosen in different ways. In \cite{SHH19}, one possible choice is the error function, $\mathtt{erf}(\cdot)$. In the numerical simulations, we choose
$$
\xi(t)=t^3\left( 10+3t\left(-5+2t\right)\right),\qquad 0\leq t\leq 1.
$$

\begin{remark} We also define an enlarged cavity $\mathcal{L}:=\Omega_0\setminus \Omega_\ses$. By the definition of $\mathcal{L}_\varepsilon$ and $\mathcal{L}_\lambda$ we have an immediate consequence of $\mathcal{L}_\lambda\subset \mathcal{L}$ and $\mathcal{L}_\varepsilon\subset \mathcal{L}$. 
\end{remark}

\begin{remark}
Here we have set the width of $\mathcal{L}_{\varepsilon}$ and $\mathcal{L}_{\lambda}$ to be $r_p$. In practice, one can use an arbitrary width (see \cite{SHH19}), but we set them to $r_p$ for simplicity.
\end{remark}

\section{Formulation of the Problem}\label{sec:form_prob}
In this section, we reduce our problem to different domains by introducing a new hybrid linear-nonlinear PB model and the domain decomposition strategy used to solve the resulting equation.

In the previous section, we introduced the dielectric permittivity and ion-exclusion functions. Let $\psi=e\tilde{\psi}/K_{\mB}T$ denote the dimensionless electrostatic potential, then the NPB equation reduces to 
\begin{equation}\label{eq:dim_npb_eqn}
-\nabla \cdot \left[\varepsilon_\abso\varepsilon(\bx)\nabla \beta\psi(\bx)\right] +\lambda(\bx)2 ce\sinh\left(\psi(\bx)\right)= \rho^{\mathrm{sol}}(\bx) \quad \mathrm{in}\ \ \mathbb{R}^3.
\end{equation}

Since the support of $\rho^{\sol}(\bx)$ is contained in $\Omega_0$, i.e., $\mathrm{supp}(\rho^{\sol})\subset \Omega_0$, by using the definition of $\varepsilon(\bx)$ and $\lambda(\bx)$, Eq.~\eqref{eq:dim_npb_eqn} reduces to
$$
-\varepsilon_s\varepsilon_\abso\beta\Delta \psi(\bx)+2 ce\sinh\left(\psi(\bx)\right) = 0\quad\mathrm{in}\ \ \Omega_{\infty}.
$$
In the bulk solvent region $\Omega_{\infty}$, we can assume that the potential $\psi$ satisfies the low potential condition, i.e.,  $|\psi|\ll 1$, and hence we can linearise the above equation to get a homogeneous screened Poisson (HSP) equation as
\begin{equation}\label{eq:hsp_eqn}
-\Delta \psi(\bx)+\kappa^2\psi(\bx)=0\quad\ \ \mathrm{in}\ \ \Omega_{\infty},
\end{equation}
where $\kappa^2=2ce/\beta\varepsilon_s\varepsilon_\abso$ is the square of the Debye H\"uckel screening constant.

Next, we note that inside the cavity $\Omega_0$, we still keep the NPB equation of the form
\begin{equation}\label{eq:nrd_eqn}
-\nabla\cdot \left[\varepsilon(\bx)\nabla \psi(\bx)\right]+\lambda(\bx)\kappa^2\varepsilon_s\sinh\left(\psi(\bx)\right)=\frac{1}{\beta \varepsilon_\abso} \rho^{\sol}(\bx)\quad\mathrm{in}\ \ \Omega_0.
\end{equation}
Along the solute-solvent boundary $\Gamma_0$, equations~\eqref{eq:hsp_eqn}--\eqref{eq:nrd_eqn} satisfy the jump conditions
\begin{equation}\label{eq:jump_cond}
\left[\![ \psi]\!\right] = 0,\qquad \left[\![\partial_{\bn}\psi]\!\right]=0\quad \mathrm{on}\quad \Gamma_0,
\end{equation}
where $\bn$ is the unit normal vector on $\Gamma_0$ pointing outward and $\partial_{\bn}=\nabla\cdot \bn$, i.e., the normal derivative which completes the hybrid linear-nonlinear model.

\begin{figure}
\centering
\begin{tikzpicture}[scale=0.75]
\fill[line width=0.8pt, color = babyblue] (1,0) -- (11,0) -- (11,8) -- (1,8) -- cycle;
\draw (1,0) -- (11,0);
\draw (11,0) -- (11,8);
\draw (11,8) -- (1,8);
\draw (1,8) -- (1,0);
\fill[line width=0.8pt, color = aliceblue] (5.5,4) circle [radius= 2.5];
\fill[line width=0.8pt, color = aliceblue] (7.5,4) circle [radius=3];
\fill[line width=0.8pt, color = white] (5.5,4) circle [radius=1];

\fill[line width=0.8pt, color = white	] (7.5,4) circle [radius=1.5];
\draw[line width=0.7pt,red] (5.875338956423, 4.926887611996) arc(247.9756872: 321.9513743:0.505);
\draw[line width=0.7pt,red] (6.42187223012, 2.957099951182) arc(44.0486257: 112.0243128:0.5);
\draw[line width=0.7pt,red] (5.875338956423, 4.926887611996) arc(67.9547342: 292.0452658:1);
\draw[line width=0.7pt,red] (6.42187223012, 2.957099951182) arc(224.0484735:494.0614503:1.5);

\begin{scriptsize}
\draw (7.5,3) node[ text width= 1.5cm, align=center] {\large{$\Omega_\ses$}};
\draw (9.7,0.4) node[ text width= 1.5cm, align=center] {\large{$\Omega_{\infty}$}};

\draw [-stealth](6,4.7) -- (6,8.2);
\draw [-stealth](9.5,4.7) -- (9,8.2);

\draw (6,8.7) node[ text width= 10cm, align=center] {\normalsize{$-\nabla\cdot \left[\varepsilon(\bx)\nabla \psi(\bx)\right]+\lambda(\bx)\kappa^2\varepsilon_s\sinh\left(\psi(\bx)\right) =\frac{1}{\beta \varepsilon_{\abso}}\rho^{\sol}(\bx)$}};

\draw (3.45,0.5) node[ text width= 4cm, align=right] {\normalsize{$-\Delta \psi(\bx) + \kappa^2\psi(\bx) = 0$}};

\draw (9.5,3.5) node[ text width= 1.5cm, align=center] {\Large{$\mathcal{L}$}};
\end{scriptsize}
\end{tikzpicture}
\caption{PDEs defined in the solute cavity $\Omega_0=\Omega_{\ses}\cup \mathcal{L}$, and the bulk solvent region $\Omega_{\infty}$.}\label{fig:pde_in_cavity}
\end{figure}
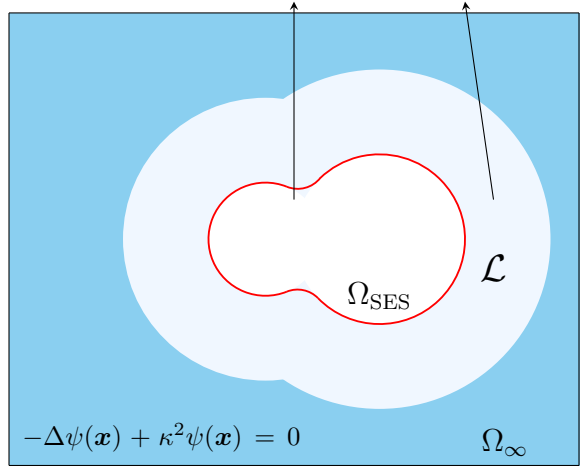

\begin{remark}
We notice that we reduced our problem defined on the whole space $\mathbb{R}^3$ to the solute cavity $\Omega_0$ and the bulk solvent region $\Omega_{\infty}$. Fig.~\ref{fig:pde_in_cavity} shows the schematic diagram of PDEs in different regions.
\end{remark}

\subsection{Transformation of the Problem}
We notice that we have two equations defined in $\mathbb{R}^3$, one inside the cavity, $\Omega_0$, and one in the bulk solvent region $\Omega_{\infty}$ in combination with interface conditions. In this subsection, we transform our problem and define them as two coupled PDEs in $\Omega_0$.

Let us denote $\psi|_{\Omega_{\infty}}$ as the Dirichlet trace from $H^1(\Omega_{\infty})$ to $L^2(\Gamma_0)$ in the sense of the trace operator, then by the linearity of Eq.~\eqref{eq:hsp_eqn} we notice that the electrostatic potential $\psi|_{\Omega_{\infty}}$ can be represented by a single-layer potential, $\widetilde{\mathcal{S}}_{\Gamma_0}:H^{-1/2}(\Gamma_0)\rightarrow H^1(\mathbb{R}^3\setminus \Gamma_0)$ as
$$
\psi(\bx)|_{\Omega_{\infty}}=\widetilde{\mathcal{S}}_{\Gamma_0}\sigma_{\me}(\bx)\quad \forall\ \bx\in\Omega_{\infty},
$$
where $\sigma_{\me}(\bx)$ is a distribution in $H^{-1/2}(\Gamma_0)$. From the continuity of the single layer operator, we can extend $\psi|_{\Omega_{\infty}}$ to $\Omega_0$ as follows:
$$
\psi_{\me}(\bx) = \int_{\Gamma_0}\frac{\exp\left(-\kappa|\bx-\by|\right)\sigma_{\me}(\by)}{4\pi |\bx-\by|}d\by\quad\forall\ \bx\in \Omega_0,
$$
where $\psi_{\me}(\bx)$ denotes the extended potential and hence $\psi_{\me}$ satisfies the HSP equation in $\Omega_{\infty}$, Eq.~\eqref{eq:hsp_eqn}, and also satisfies
$$
-\Delta\psi_{\me}(\bx)+\kappa^2\psi_{\me}(\bx)=0\quad \mathrm{in}\ \ \Omega_0.
$$
Here we also introduce the single-layer potential operator $\mathcal{S}_{\Gamma_0}:H^{-1/2}(\Gamma_0)\rightarrow H^{1/2}(\Gamma_0)$ defined by
\begin{equation}\label{eq:single_layer_invert}
\mathcal{S}_{\Gamma_0}\sigma_{\me}(\bx):=\int_{\Gamma_0}\frac{\exp\left(-\kappa|\bx-\by|\right)\sigma_{\me}(\by)}{4\pi|\bx-\by|}d\by\quad\forall \bx\in \Gamma_0,
\end{equation}
which is an invertible operator and hence $\sigma_{\me}=\mathcal{S}_{\Gamma_0}^{-1}\psi|_{\Gamma_0}$.

From \cite[Theorem~3.3.1]{SS11} we also have the relationship
$$
\sigma_{\me} = \partial_{\bn}\psi_{\me}|_{\Omega_0}-\partial_{\bn}\psi_{\me}|_{\Omega_{\infty}}\quad \mathrm{on}\quad\Gamma_0.
$$
By the jump condition (Eq.~\eqref{eq:jump_cond}) of the normal derivative we have
$$
\partial_{\bn}\psi|_{\Omega_{\infty}}-\partial_{\bn}\psi|_{\Omega_0}=0\quad\ \ \mathrm{on}\ \ \Gamma_0,
$$
which implies $\partial_{\bn}\psi_{\me}|_{\Omega_{\infty}}=\partial_{\bn}\psi|_{\Omega_0}$ on $\Gamma_0$.  Hence,
$$
\sigma_{\me}=\partial_{\bn}\left(\psi_{\me}|_{\Omega_0}-\psi|_{\Omega_0}\right)\quad \mathrm{on}\ \ \Gamma_0.
$$
Also, by the jump condition of the potential, we get $\psi_{\me}=\psi|_{\Omega_0}$ on $\Gamma_0$. Hence, the extended potential $\psi_{\me}$ is defined on $\Omega_0$ by
\begin{alignat*}{3}
-\Delta \psi_{\me}(\bx)+\kappa^2\psi_{\me}(\bx)&=0&&\quad \mathrm{in}\quad\Omega_0,\nonumber\\
\psi_{\me}(\bx)& = \psi(\bx) &&\quad\mathrm{on}\quad \Gamma_0.
\end{alignat*}

Now, we move towards Eq.~\eqref{eq:nrd_eqn}. Let $\psi_0(\bx)$ denote the potential generated by $\rho^{\sol}(\bx)/\beta \varepsilon_\abso$ in vacuum (see Remark~\ref{rem:rho_sol_def}), i.e.,
$$
\psi_0(\bx) = \sum_{i=1}^M\frac{q_i}{4\pi \varepsilon_\abso\beta|\bx-\bx_i|},
$$
satisfying
\begin{equation}\label{eq:potential_vacuum}
-\Delta\psi_0=\frac{1}{\beta \varepsilon_\abso} \rho^{\sol}(\bx)\quad\mathrm{in}\ \ \mathbb{R}^3.
\end{equation}

Let us denote the reaction potential by $\psi_{\mr}:=\psi-\psi_0$, i.e., the difference between the electrostatic potential with and without the presence of a solvent. Then Eq.~\eqref{eq:nrd_eqn} equivalently writes
\begin{alignat}{3}\label{eq:temp_e_1}
-\nabla\cdot \left[\varepsilon(\bx)\nabla \psi_{\mr}(\bx)\right] &+ \lambda(\bx)\kappa^2\varepsilon_s\sinh\left(\left(\psi_{\mr}+\psi_0\right)(\bx)\right) \nonumber \\
&= \frac{1}{\beta\varepsilon_\abso} \rho^{\sol}(\bx)+\nabla\cdot\left[\varepsilon(\bx)\nabla\psi_0(\bx)\right]&&\qquad\mathrm{in}\ \ \Omega_0.
\end{alignat}
Substituting Eq.~\eqref{eq:potential_vacuum} in Eq.~\eqref{eq:temp_e_1} and denoting $\sinh(\Phi)$ by $\mathcal{F}(\Phi)\Phi$ where
$$
\mathcal{F}(\Phi)=\frac{\sinh(\Phi)}{\Phi},
$$
for any positive function $\Phi$. We further reduce the equation to
\begin{alignat*}{3}
-\nabla\cdot \left[\varepsilon(\bx)\nabla \psi_{\mr}(\bx)\right]&+\lambda(\bx)\kappa^2\varepsilon_s\mathcal{F}\left( (\psi_{\mr}+\psi_0)(\bx)\right)\left(\psi_{\mr}+\psi_0\right)(\bx)\nonumber \nonumber\\
&=\nabla\cdot\left[\left(\varepsilon(\bx)-1\right)\nabla\psi_0(\bx)\right]\qquad&&\mathrm{in}\ \ \Omega_0.
\end{alignat*}

From the jump condition, $\left[\![\psi]\!\right]=0$ we have $\psi|_{\Omega_0}-\psi|_{\Omega_{\infty}}=0$ which implies
$$
\psi_{\mr}+\psi_0=\psi_{\me}\quad\mathrm{on}\ \ \Gamma_0,
$$
and similarly
$$
\sigma_{\me}=\partial_{\bn}\psi_{\me}-\partial_{\bn}\left(\psi_{\mr}+\psi_0\right)\quad\mathrm{on}\ \ \Gamma_0.
$$
From Eq.~\eqref{eq:single_layer_invert}, we also get the global coupling condition
$$
\psi_{\me}(\bx)=\mathcal{S}_{\Gamma_0}\sigma_{\me}=\mathcal{S}_{\Gamma_0}\left[ \partial_{\bn}\psi_{\me}-\partial_{\bn}\left(\psi_{\mr}+\psi_0\right)\right]\quad\mathrm{on}\ \ \Gamma_0.
$$

In summary, our original equation reduces to
\begin{alignat}{3}
-\nabla\cdot \left[\ \varepsilon(\bx)\nabla \psi_{\mr}(\bx)\right]&+\lambda(\bx)\kappa^2\varepsilon_s\mathcal{F}\left((\psi_{\mr}+\psi_0)(\bx)\right)\left(\psi_{\mr}+\psi_0\right)(\bx) \nonumber \\ &=\nabla\cdot\left[\left(\varepsilon(\bx)-1\right)\nabla\psi_0(\bx)\right]&&\mathrm{in}\ \ \Omega_0,\label{eq:coupled_pde_1}\\
-\Delta \psi_{\me}(\bx)+\kappa^2\psi_{\me}(\bx) &= 0 &&\mathrm{in}\ \ \Omega_0,\label{eq:coupled_pde_2}
\end{alignat}
with two coupling conditions given on $\Gamma_0$ by
\begin{alignat}{3}
\left(\psi_{\mr}+\psi_0\right)(\bx)  &=  \psi_{\me}(\bx) &&\quad \mathrm{on}\ \ \Gamma_0,\label{eq:coupled_bc_1}\\
\psi_{\me}(\bx) &=  g &&\quad \mathrm{on}\ \ \Gamma_0,\label{eq:coupled_bc_2}
\end{alignat}
where
\begin{equation}\label{eq:global_coupling}
g  =  \mathcal{S}_{\Gamma_0}\sigma_{\me}  = \mathcal{S}_{\Gamma_0}\left[ \partial_{\bn}\psi_{\me}-\partial_{\bn}\left(\psi_{\mr}+\psi_0\right)\right]\quad\mathrm{on}\ \ \Gamma_0.
\end{equation}

\begin{remark}\label{rem:limiting_case}
The PDEs defined in Eq.~\eqref{eq:coupled_pde_1}--\eqref{eq:coupled_pde_2} encompass different solvation models. As mentioned in the introduction, we can linearize Eq.~\eqref{eq:coupled_pde_1} using $\mathcal{F}(\Phi)\approx 1$, if $|\Phi|\ll 1$ and get the LPB equation based on the SES. In \cite{QSM19}, a domain-decomposition method was developed for this model but on a vdW-cavity (called ddLPB). If $\kappa\rightarrow 0$, we recover the polarizable continuum model with SES boundary (PCM-SES). We get the classical PCM model for the vdW-cavity ($r_p\rightarrow 0$, $a\rightarrow 0$). Lastly, for $\kappa\rightarrow \infty$, the model tends to the conductor line screening model (COSMO), which is reasonable as then the solvent becomes a perfect conductor as the ionic strength tends to $\infty$ and screens any change from the solute. Domain decomposition algorithms for all the three models denoted by ddPCM-SES, ddPCM, and ddCOSMO can be found in \cite{QSM18,  SCLM16, CMS13}, respectively. 
\end{remark}

\subsection{Global Strategy}\label{sec:global_strategy}

For solving Eq.~\eqref{eq:coupled_pde_1}--\eqref{eq:coupled_pde_2}, we follow the classical iterative strategy developed for solvation models as described in \cite{QSM18, QSM19}. Let $g^{(0)}$ be an initial guess for the Dirichlet condition $\psi_{\me}|_{\Gamma_0}$ and set $\mk=1$:
\begin{enumerate}
\item[$\mathtt{[1]}$:] Solve the following nonlinear Dirichlet boundary problem for $\psi_{\mr}^{(\mk)}$:
\begin{alignat*}{3}
-\nabla\cdot \left[ \varepsilon(\bx)\nabla \psi_{\mr}^{(\mk)}(\bx)\right]&+\lambda(\bx)\kappa^2\varepsilon_s\mathcal{F}\left(\psi_{\mr}^{(\mk)}+\psi_0\right) \left(\psi_{\mr}^{(\mk)}+\psi_0\right)(\bx) \\
& = \nabla\cdot\left[\left(\varepsilon(\bx)-1\right)\nabla\psi_0(\bx)\right]  \qquad &&\mathrm{in}\ \ \Omega_0,\\
\psi_{\mr}^{(\mk)}(\bx) &= g^{(\mk-1)} -\psi_0(\bx) \qquad &&\mathrm{on}\ \  \Gamma_0,
\end{alignat*}
and derive the Neumann trace $\partial_{\bn}\psi_{\mr}^{(\mk)}$ on $\Gamma_0$.
\item[$\mathtt{[2]}:$] Solve the Dirichlet boundary problem for $\psi_{\me}^{(\mk)}$:
\begin{alignat*}{3}
-\Delta \psi_{\me}^{(\mk)}(\bx)+\kappa^2\psi_{\me}^{(\mk)}(\bx)&  = 0 \qquad && \mathrm{in}\ \ \Omega_0,\\
\psi_{\me}^{(\mk)}(\bx) &= g^{(\mk-1)} \qquad &&\mathrm{on}\ \ \Gamma_0,
\end{alignat*}
and derive the Neumann trace $\partial_{\bn}\psi_{\me}^{(\mk)}$ on $\Gamma_0$.
\item[$\mathtt{[3]}:$] Build the charge density $\sigma_{\me}^{(\mk)}=\partial_{\bn}\psi_{\me}^{(\mk)}-\partial_{\bn}\left(\psi_{\mr}^{(\mk)}+\psi_0\right)$ and compute a new Dirichlet condition $g^{(\mk)}=\mathcal{S}_{\Gamma_0}\sigma_{\me}^{(\mk)}$.
\item[$\mathtt{[4]}:$] Compute the contribution $E^{\mk}_s$ to the solvation energy based on $\psi_{\mr}^{(\mk)}$ at the $\mk^{\ith}$ iteration step, set $\mk\rightarrow \mk+1$, go back to $\mathtt{[1]}$ and repeat until $|E^{\mk}_s-E^{\mk-1}_s|/|E^{\mk}_s|<\mathtt{tol}$ for $\mathtt{tol}\ll 1$.
\end{enumerate}
In the above algorithm, $E^{\mk}_s$ denotes the electrostatic solvation energy at iteration $\mk$, which will be defined in Sec.~\ref{sec:numres}.

\begin{remark} For choosing a suitable guess to $g^{(0)}$ (defined on $\Gamma_0$), we consider the (unrealistic) situation when the whole space $\mathbb{R}^3$ is covered by the solvent medium. Then, the electrostatic potential $\psi$ would be given by
$$
\psi(\bx)=\sum_{i=1}^M\frac{q_i}{\beta\varepsilon_s}\frac{\exp\left(-\kappa |\bx-\bx_i|\right)}{4\pi \varepsilon_\abso|\bx-\bx_i|}\quad \forall\ \ \bx\in \mathbb{R}^3,
$$
see \cite[Sec~1.3.2]{Lamm03}. Hence $g^{(0)}$ is chosen as this potential restricted to $\Gamma_0$.
\end{remark}

\begin{remark} We note that the global strategy is an iterative process. After discretization, the final convergent solution satisfies a global nonlinear system that can be solved. Unlike the LPB approach, the solution of the nonlinear problem in Step~$\mathtt{[1]}$ needs to be studied properly and will be discussed in Remark~\ref{rem:linear}.
\end{remark}

\subsection{Domain Decomposition Strategy}

In this work, we will consider the Schwarz domain decomposition method \cite{QV99} as it aims to solve PDEs defined on complex domains, which can be decomposed as a union of overlapping and simple sub-domains. The main idea is to solve the same equation in each sub-domain but with boundary conditions that depend on the global boundary conditions and the solution of neighboring domains.

We recall that we have a natural decomposition of $\Omega_0$ as follows
$$
\Omega_0=\bigcup_{j=1}^M\Omega_j,\quad \Omega_j=B_{R_j}(\bx_j).
$$

\begin{figure}[t!]
\begin{center}
\begin{tikzpicture}[scale=0.75]
\draw[line width=0.8pt] (1,1) circle [radius=1];
\draw[line width=0.8pt] (2.5,1) circle [radius=1.5];
\draw[line width=0.8pt] (4,1.5) circle [radius=1];
\draw[line width=0.8pt] (4,0.5) circle [radius=1];
\draw[line width=0.8pt,red] (3.3345262,2.24642125) arc(56.2:141.06:1.5);
\draw[line width=0.8pt,red] (1.333334, 0.05719) arc(218.94:303.8:1.5);
\draw[line width=0.8pt,blue] (1.333334, 1.942809) arc(141.06:218.94:1.5);
\draw[line width=0.8pt,blue] (3.334526249,-0.246421253) arc(303.8:416.2:1.5);
\begin{scriptsize}
\draw (2.5,1) node[align=center] {\small{$\Omega_j$}};
\draw (2.5,1) node[align=center] {\small{$\Omega_j$}};
\draw (2,2.7) node[align=center] {\small{$\textcolor{red}{\Gamma_j^{\te}} = \Gamma_0\cap \Gamma_j$}};
\draw (5.1,0) node[align=right] {\small{$\textcolor{blue}{\Gamma_j^{\ti}} = \Omega_0\cap \Gamma_j$}};
\end{scriptsize}
\end{tikzpicture}
\caption{2-D schematic diagram of $\Gamma_j^{\ti}$ and $\Gamma_j^{\te}$.}\label{fig:external_internal}
\end{center}
\end{figure}
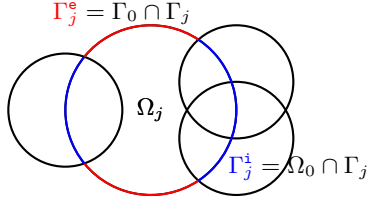

We replace the global equation~\eqref{eq:coupled_pde_1}--\eqref{eq:coupled_pde_2} by the following coupled equations, each restricted to $\Omega_j$,
\begin{alignat}{3}\label{eq:ndr_subdomain}
-\nabla\cdot \left[ \varepsilon(\bx)\nabla \psi_{\mr}|_{\Omega_j}(\bx)\right]&+\lambda(\bx)\kappa^2\varepsilon_s\mathcal{F}\left(\psi_{\mr}|_{\Omega_j}+\psi_0\right) \left(\psi_{\mr}|_{\Omega_j}+\psi_0\right)(\bx)\nonumber\\
& = \nabla\cdot\left[\left(\varepsilon(\bx)-1\right)\nabla\psi_0(\bx)\right] \qquad &&\mathrm{in}\ \ \Omega_j,\nonumber\\
\psi_{\mr}|_{\Omega_j}(\bx) &=h_{\mr,j}\quad &&\mathrm{on}\ \ \Gamma_j,
\end{alignat}
with
\begin{equation}\label{eq:ndr_bc}
h_{\mr,j}=\left\lbrace
\begin{array}{lll}
\psi_{\mr} & \mathrm{on} & \Gamma_j^{\ti},\\
g-\psi_0 & \mathrm{on} & \Gamma_j^{\te}.
\end{array}\right.
\end{equation}
Here, $\Gamma_j^{\te}$ is the external part of $\Gamma_j$ not contained in any other $\Omega_i\ (i\neq j)$, i.e., $\Gamma_j^{\te}=\Gamma_0\cap \Gamma_j$, and $\Gamma_j^{\ti}$ is the internal part of $\Gamma_j$, i.e., $\Gamma_j^{\ti}=\Omega_0\cap \Gamma_j$ (see Fig.~\ref{fig:external_internal}). 

Similarly, the equations for the extended potential are
\begin{alignat}{3}\label{eq:hsp_subdomain}
-\Delta \psi_{\me}|_{\Omega_j}(\bx)+\kappa^2\psi_{\me}|_{\Omega_j}(\bx)&  = 0\qquad &&\mathrm{in}\ \ \Omega_j,\nonumber\\
\psi_{\me}|_{\Omega_j}(\bx) &= h_{\me,j} \qquad &&\mathrm{on}\ \ \Gamma_j,
\end{alignat}
with
\begin{equation}\label{eq:hsp_bc}
h_{\me,j}=\left\lbrace
\begin{array}{lll}
\psi_{\me} & \mathrm{on} & \Gamma_j^{\ti},\\
g & \mathrm{on} & \Gamma_j^{\te}.
\end{array}\right.
\end{equation}

\begin{figure}[t!]
\centering
\begin{tikzpicture}[scale=1.50]
\fill[line width=0.8pt, color = white] (0,0) -- (8,0) -- (8,7) -- (0,7) -- cycle;
\fill[line width=0.8pt] (0.6,1) circle [radius=0.35];
\fill[line width=0.8pt] (1.5,1) circle [radius=0.45];
\fill[line width=0.8pt] (2.25,1.35) circle [radius=0.3];
\fill[line width=0.8pt] (2.25,0.65) circle [radius=0.3];
\fill[line width=0.8pt, color = white] (0.6,1) circle [radius=0.34];
\fill[line width=0.8pt, color = white] (1.5,1) circle [radius=0.44];
\fill[line width=0.8pt, color = white] (2.25,1.35) circle [radius=0.29];
\fill[line width=0.8pt, color = white] (2.25,0.65) circle [radius=0.29];
\fill[line width=0.8pt] (1,3.5) circle [radius=0.35];
\fill[line width=0.8pt] (1.5,3.5) circle [radius=0.45];
\fill[line width=0.8pt] (1.85,3.85) circle [radius=0.3];
\fill[line width=0.8pt] (1.85,3.15) circle [radius=0.3];
\fill[line width=0.8pt, color = white] (1,3.5) circle [radius=0.34];
\fill[line width=0.8pt, color = white] (1.5,3.5) circle [radius=0.44];
\fill[line width=0.8pt, color = white] (1.85,3.85) circle [radius=0.29];
\fill[line width=0.8pt, color = white] (1.85,3.15) circle [radius=0.29];
\fill[line width=0.8pt] (5.8,1) circle [radius=0.35];
\fill[line width=0.8pt] (6.7,1) circle [radius=0.45];
\fill[line width=0.8pt] (7.45,1.35) circle [radius=0.3];
\fill[line width=0.8pt] (7.45,0.65) circle [radius=0.3];
\fill[line width=0.8pt, color = babyblue] (5.8,1) circle [radius=0.34];
\fill[line width=0.8pt, color = babyblue] (6.7,1) circle [radius=0.44];
\fill[line width=0.8pt, color = babyblue] (7.45,1.35) circle [radius=0.29];
\fill[line width=0.8pt, color = babyblue] (7.45,0.65) circle [radius=0.29];
\fill[line width=0.8pt] (6.2,3.5) circle [radius=0.35];
\fill[line width=0.8pt] (6.7,3.5) circle [radius=0.45];
\fill[line width=0.8pt] (7.05,3.85) circle [radius=0.3];
\fill[line width=0.8pt] (7.05,3.15) circle [radius=0.3];
\fill[line width=0.8pt, color = babyblue] (6.2,3.5) circle [radius=0.34];
\fill[line width=0.8pt, color = babyblue] (6.7,3.5) circle [radius=0.44];
\fill[line width=0.8pt, color = babyblue] (7.05,3.85) circle [radius=0.29];
\fill[line width=0.8pt, color = babyblue] (7.05,3.15) circle [radius=0.29];
\fill[line width=0.8pt, color = babyblue] (2.8,5.05) -- (5,5.05) -- (5,6.95) -- (2.8,6.95) -- cycle;
\fill[line width=0.8pt] (3.5,6) circle [radius=0.35];
\fill[line width=0.8pt] (4,6) circle [radius=0.45];
\fill[line width=0.8pt] (4.35,6.35) circle [radius=0.3];
\fill[line width=0.8pt] (4.35,5.65) circle [radius=0.3];
\fill[line width=0.8pt, color = white] (3.5,6) circle [radius=0.34];
\fill[line width=0.8pt, color = white] (4,6) circle [radius=0.44];
\fill[line width=0.8pt, color = white] (4.35,6.35) circle [radius=0.29];
\fill[line width=0.8pt, color = white] (4.35,5.65) circle [radius=0.29];
\draw[line width=0.8pt,red] (1.2,1.34) arc(131.81:228.19:0.45);
\draw[line width=0.8pt,red] (0.77,0.69) arc(299.06:420.94:0.35);

\draw[line width=0.8pt,red] (1.9452,1.0655) arc(8.367:87.89:0.45);
\draw[line width=0.8pt,red] (1.56551,0.55479) arc(278.37:351.633:0.45);

\draw[line width=0.8pt,red] (1.96551,1.44521) arc(161.5:251.5:0.3);
\draw[line width=0.8pt,red] (2.34521,0.93449) arc(71.5:199.5:0.3);
\draw[line width=0.8pt,red] (6.4,1.34) arc(131.81:228.19:0.45);
\draw[line width=0.8pt,red] (5.97,0.69) arc(299.06:420.94:0.35);

\draw[line width=0.8pt,red] (7.1452,1.0655) arc(8.367:87.89:0.45);
\draw[line width=0.8pt,red] (6.76551,0.55479) arc(278.37:351.633:0.45);

\draw[line width=0.8pt,red] (7.16551,1.44521) arc(161.5:251.5:0.3);
\draw[line width=0.8pt,red] (7.54521,0.93449) arc(71.5:199.5:0.3);
\fill[line width=0.8pt, color = aliceblue] (2.7,3.2) -- (5.3,3.2) -- (5.3,3.8) -- (2.7,3.8) -- cycle;
\draw (2.15,3.5) -- (2.65,3.5);
\draw (5.35,3.5) -- (5.8,3.5);
\draw [-stealth](2.7,6)-- (1.5,4.5);
\draw [-stealth](5.1,6)-- (6.5,4.5);
\draw [-stealth](1.5,2.75)-- (1.5,1.7);
\draw [-stealth](6.5,2.75)-- (6.5,1.7);
\begin{scriptsize}
\draw (3.6,6) node[ text width= 1.5cm, align=right] {Eq.~\eqref{eq:nrd_eqn}};
\draw (3.8,6.8) node[ text width= 2.95cm, align=left] {Eq.~\eqref{eq:hsp_eqn}};
\draw (3.3,5.5) node[ text width= 1.5cm, align=right] {$\left[\![\psi]\!\right]  =  0$};
\draw (3.3,5.2) node[ text width= 1.5cm, align=right] {$\left[\![\partial_{\bn}\psi]\!\right]  =  0$};
\draw (4.45,5.2) node[ text width= 1.5cm, align=right] {$\Omega_{\infty}$};
\draw (4.0,6.5) node[ text width= 1.5cm, align=right] {$\Omega_0$};
\draw (4.2,6) node[ text width= 1.5cm, align=right] {$\Gamma_0$};
\draw (1.1,3.5) node[ text width= 2.5cm, align=right] {Eq.~\eqref{eq:coupled_pde_1}--\eqref{eq:coupled_bc_1}};
\draw (6.3,3.5) node[ text width=2.5cm, align=right] {Eq.~\eqref{eq:coupled_pde_2}--\eqref{eq:coupled_bc_2}};
\draw (3.97,3.5) node[text width = 4cm, align = right] {\scriptsize{$
g  = \mathcal{S}_{\Gamma_0}\left[ \partial_{\bn}\psi_e-\partial_{\bn}\left(\psi_r+\psi_0\right)\right]
$}};
\end{scriptsize}
\end{tikzpicture}
\caption{Schematic diagram of the domain decomposition algorithm for the NPB equation.}\label{fig:ddpb_model}
\end{figure}
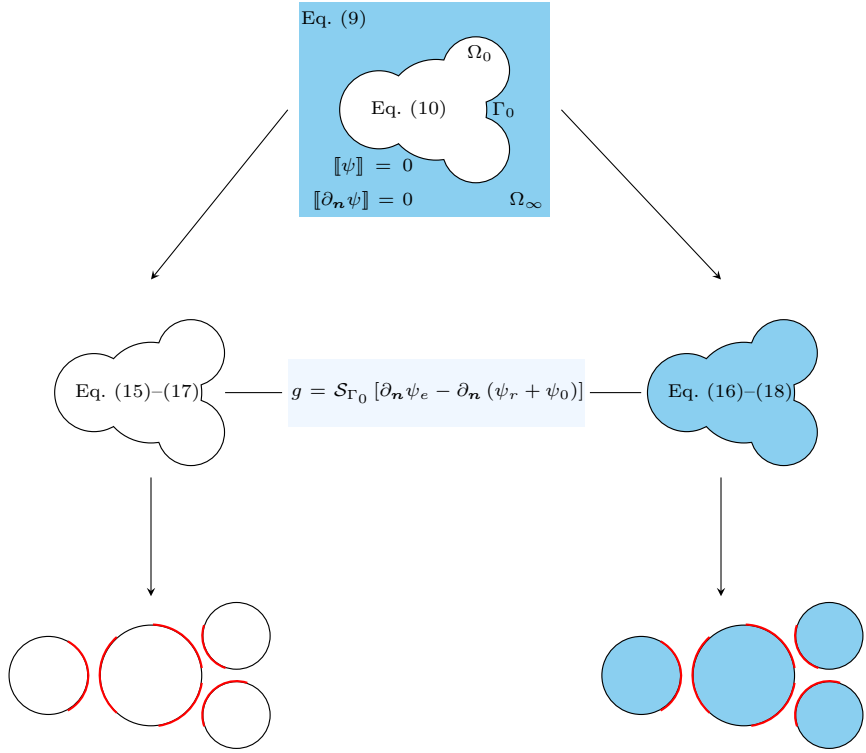

\begin{remark} The Dirichlet boundary conditions in Eq.~\eqref{eq:ndr_subdomain}--\eqref{eq:hsp_subdomain} are implicit since $\psi_{\mr}$ (respectively $\psi_{\me}$) is not known on $\Gamma_j^{\ti}$. Hence, one needs to use an iterative scheme to solve Eq.~\eqref{eq:ndr_subdomain}--\eqref{eq:ndr_bc} (respectively Eq.~\eqref{eq:hsp_subdomain}--\eqref{eq:hsp_bc}), such as the parallel Schwarz algorithm or the alternating Schwarz algorithm as presented in for ddCOSMO in \cite{CMS13}. The visual overview of the domain decomposition algorithm is presented in Fig.~\ref{fig:ddpb_model}.
\end{remark}

\begin{remark}\label{rem:linear}
We note that in Eq.~\eqref{eq:ndr_subdomain} we have a nonlinearity because of the term
\begin{equation}\label{eq:nonlinear_term}
\lambda(\bx)\kappa^2\varepsilon_s\mathcal{F}\left(\psi_{\mr}|_{\Omega_j}+\psi_0\right)\left(\psi_{\mr}|_{\Omega_j}+\psi_0\right)(\bx).
\end{equation}
A standard way of solving such nonlinear equations after discretization is the use of a fixed point technique and replacing Eq.~\eqref{eq:nonlinear_term} by
$$
\lambda(\bx)\kappa^2\varepsilon_s\mathcal{F}\left(\psi_{\mr}|_{\Omega_j}^{(\nu-1)}+\psi_0\right)\left(\psi_{\mr}|_{\Omega_j}^{(\nu)}+\psi_0\right)(\bx).
$$
where $\psi_{\mr}^{(\nu-1)}(\bx)$ denotes the solution at the $(\nu-1)^{\ith}$ iterative step.  As $\psi_{\mr}^{(\nu-1)}(\bx)$ is known at the $(\nu)^{\ith}$ iterative step, we can replace Eq.~\eqref{eq:ndr_subdomain} by the linear counterpart
\begin{alignat}{3}\label{eq:ndr_subdomain_linear}
-\nabla\cdot \left[ \varepsilon(\bx)\nabla \psi_{\mr}|_{\Omega_j}(\bx)\right]&+\lambda(\bx)\kappa^2\varepsilon_s\mathcal{F}\left((\overline{\psi_{\mr}}|_{\Omega_j}+\psi_0)\right) \left(\psi_{\mr}|_{\Omega_j}+\psi_0\right)(\bx)\nonumber\\
& = \nabla\cdot\left[\left(\varepsilon(\bx)-1\right)\nabla\psi_0(\bx)\right]\qquad &&\mathrm{in}\ \ \Omega_j,\nonumber\\
\psi_{\mr}|_{\Omega_j}(\bx) &=h_{\mr,j}\qquad &&\mathrm{on}\ \ \Gamma_j,
\end{alignat}
where we have dropped the notation of $(\nu)^{\ith}$ iterative step and denote $(\nu-1)^{\ith}$ solution of $\psi_{\mr}(\bx)$ by $\overline{\psi_{\mr}}(\bx)$. In the subsequent sections, we refer to Eq.~\eqref{eq:ndr_subdomain_linear} as a generalized screened Poisson (GSP) equation.
\end{remark}

\section{Single Domain Solvers}\label{sec:single_domain}
In this section, we will develop two single domain solvers in the unit ball for solving Eq.~\eqref{eq:hsp_subdomain} and Eq.~\eqref{eq:ndr_subdomain_linear}, respectively. Without loss of generality, we will consider a unit ball centered at the origin for developing the solvers.

\subsection{HSP Solver}
In \cite{QSM19}, an HSP solver was developed for the ddLPB meth-od. One can use the same solver for ddNPB as well. For completeness, we outline the main ideas in this subsection.

The HSP equation in the unit ball is given by
\begin{alignat}{3}\label{eq:hsp_in_unit}
-\Delta u_{\me} + \kappa^2u_{\me} & = 0 \qquad&& \mathrm{in}\ \ \ B_1(\bo),\nonumber \\
u_{\me} & = \phi_{\me} \qquad && \mathrm{on}\ \ \ \ssphere,
\end{alignat}
where $B_1(\bo)$ denotes the unit ball at center $\bo$ and for $j=1,\dots, M$, $\phi_{\me}(\bx) = h_{\me,j}(\bx_j+R_j \bx)$ for the HSP equation in the sub-domain $\Omega_j$.

The solution of Eq.~\eqref{eq:hsp_in_unit} in $H^1(B_1(\bo))$ can be written as
$$
u_{\me}(r,\theta, \varphi) = \sum_{\ell =  0}^{\infty}\sum_{m=-\ell}^{\ell}\left[ \phi_{\me}\right]_{\ell}^m\frac{i_{\ell}(r)}{i_{\ell}(1)}Y_{\ell}^m(\theta,\varphi),\quad 0\leq r\leq 1,\ \  0\leq \theta\leq \pi,\quad 0\leq \varphi \leq 2\pi,
$$
where $i_{\ell}$ is the modified spherical Bessel function of the first kind, see \cite[Chapter~14]{AWH13}, $Y_{\ell}^m$\ denotes the (real orthonormal) spherical harmonics of degree $\ell$ and order $m$ defined on $\ssphere$, and
$$
\left[\phi_{\me}\right]_{\ell}^m=\int_{\partial B_1(\bo)}\phi_{\me}(\bs)Y_{\ell}^m(\bs)d\bs,
$$
is the real coefficient of $u_{\me}$ corresponding to the mode $Y_{\ell}^m$.  Now $u_{\me}$ can be numerically approximated by $\overline{u}_{\me}$ in the discretization space spanned by truncated basis of spherical harmonics $\left\lbrace Y_{\ell}^m\right\rbrace_{0\leq \ell\leq \ell_{\max},-\ell\leq m\leq \ell}$, defined by
\begin{equation}\label{eq:hsp_sol_unit}
\overline{u}_{\me}(r,\theta,\varphi)=\sum_{\ell, m}\left[ \tilde{\phi}_{\me}\right]_{\ell}^m\dfrac{i_{\ell}(r)}{i_{\ell}(1)}Y_{\ell}^m(\theta,\varphi),\quad 0\leq r\leq 1,\quad 0\leq \theta\leq \pi,\quad 0\leq \varphi \leq 2\pi,
\end{equation}
where
$$
\sum_{\ell, m}=\sum_{\ell=0}^{\ell_{\max}}\sum_{m=-\ell}^{\ell},
$$
$\ell_{\max}$ denotes the maximum degree of spherical harmonics and
$$
\left[\tilde{\phi}_{\me}\right]_{\ell}^m=\sum_{n=1}^{N_{\leb}}\omega_n^{\leb}\phi_{\me}(\bs_n)Y_{\ell}^m(\bs_n).
$$
To approximate the integration, we use the Lebedev quadrature formula \cite{Hax07} where $\bs_n\in \ssphere$ are the Lebedev quadrature points \cite{Hax07}, $\omega_n^{\leb}$ are the corresponding weights, and $N_{\leb}$ are the number of Lebedev quadrature points. The reason we use Lebedev quadrature over other quadrature is that it gives exact integration for spherical harmonics up to a given order (depending on the quadrature points), achieving thus good accuracy for few point evaluations, see \cite{CMS13}.

\subsection{GSP Solver}\label{sec:gsp_solver}
For the GSP equation~\eqref{eq:ndr_subdomain_linear}, consider the following problem in the unit ball
\begin{alignat}{3}\label{eq:ndr_unit}
-\nabla\cdot\left[ \tilde{\varepsilon}(\bx)\nabla u(\bx)\right]+\tilde{\lambda}(\bx)\widetilde{\mathcal{F}}\left(\overline{u}(\bx)\right)u(\bx) & = f(\bx)\qquad && \mathrm{in}\ \ \ B_1(\bo),\nonumber\\
u(\bx) & = \phi_{\mr}(\bx)\qquad && \mathrm{on}\ \ \ \ssphere,
\end{alignat}
where for $j=1,\dots, M$,  $\widetilde{\mathcal{F}}\left(\overline{u}(\bx)\right)=\kappa^2\varepsilon_s\mathcal{F}\left(\left(\overline{\psi_{\mr}}+\psi_0\right)\left(\bx_j+R_j\bx\right)\right)$; $f(\bx)=\nabla\cdot\left[\left(\varepsilon(\bx_j+R_j\bx)-1\right)\nabla \psi_0\left(\bx_j+R_j\bx\right)\right]-\mathcal{F}\left(\left(\overline{\psi_{\mr}}+\psi_0\right)\left(\bx_j+R_j\bx\right)\right)\psi_0(\bx_j+R_j\bx)$; $\tilde{\varepsilon}(\bx)\linebreak=\varepsilon(\bx_j+ R_j\bx)$,  $\tilde{\lambda}(\bx)=\lambda(\bx_j+ R_j\bx)$, and $\phi_{\mr}(\bx) = h_{\mr,j}\left(\bx_j+R_j\bx\right)$ for the GSP equation in the sub-domain $\Omega_j$.

For a Laplace equation (see \cite{QSM18}) it is known that there exists a $\hat{u}_1\in H^1\left(B_1(\bo)\right)$ such that
\begin{alignat}{3}\label{eq:harmonic_eq_1}
-\Delta \hat{u}_1 & = 0\quad &&\mathrm{in}\quad B_1(\bo),\nonumber\\
\hat{u}_1 & = \phi_{\mr}\quad && \mathrm{on}\quad\ssphere,
\end{alignat}
and $\hat{u}_1$ can be approximated by spherical harmonics in a similar way as Eq.~\eqref{eq:hsp_sol_unit} (see \cite{CMS13} for the exact representation). Let $w=u-\hat{u}_1\in H_0^1\left( B_1(\bo)\right)$ and as a consequence $w$ satisfies
\begin{alignat*}{3}
-\nabla \cdot \left[\tilde{\varepsilon}(\bx)\nabla\left(w+\hat{u}_1\right)(\bx)\right]+\tilde{\lambda}(\bx)\widetilde{\mathcal{F}}\left(\left(\overline{w+\hat{u}_1}\right)\right)\left(w+\hat{u}_1\right)(\bx) &= f(\bx)\ \ &&\mathrm{in}\ \ B_1(\bo),\nonumber\\
w(\bx) & =  0\ \  &&\mathrm{on}\ \ \ssphere,
\end{alignat*}
which can be further reduced to
\begin{alignat}{3}\label{eq:harmonic_eq_2}
-\nabla \cdot \left[\tilde{\varepsilon}(\bx)\nabla w(\bx)\right]+\tilde{\lambda}(\bx)\widetilde{\mathcal{F}}\left(\left(\overline{w+\hat{u}_1}\right)\right)w(\bx) & = \tilde{f}(\bx)\quad &&\mathrm{in}\ \ B_1(\bo),\nonumber\\
w(\bx) & =  0\quad &&\mathrm{on}\ \ \ssphere,
\end{alignat}
where $\tilde{f}(\bx) = f(\bx) + \nabla\cdot \left[\tilde{\varepsilon}(\bx)\nabla \hat{u}_1(\bx)\right] - \tilde{\lambda}(\bx)\widetilde{\mathcal{F}}\left(\left(\overline{w+\hat{u}_1}\right)(\bx)\right)\hat{u}_1(\bx)$.  Throughout this section we use $\widetilde{\mathcal{F}}\left(\overline{w}^{\hat{u}_1}(\bx)\right)$ to denote $\widetilde{\mathcal{F}}\left(\left(\overline{w+\hat{u}_1}\right)(\bx)\right)$.

From the definition of $\varepsilon(\bx), \tilde{\lambda}(\bx)$, and the fact that the vdW-ball, $B_{r_j}(\bx_j)\subset \Omega_j$, the reaction potential $\psi_{\mr}$ defined in Eq.~\eqref{eq:ndr_subdomain} is harmonic in the smaller ball $B_{r_j}(\bx_j)$ and consequently $u(\bx)$ in Eq.~\eqref{eq:ndr_unit} and $w(\bx)$ in Eq.~\eqref{eq:harmonic_eq_2} are harmonic in $B_{\delta}(\bo)$, where
$$
\delta=\frac{r_j}{r_j+r_p+r_0+a}\in (0,1).
$$
Let $\mathcal{D}_{\delta}=B_1(\bo)\setminus B_{\delta}(\bo)$ be the region between $\ssphere$ and $\partial B_{\delta}(\bo)$ and we define the subspace $H_{0, \delta}^1(\mathcal{D}_{\delta})\subset H^1(\mathcal{D}_{\delta})$ as follows
$$
H_{0, \delta}^1(\mathcal{D}_{\delta})=\left\lbrace w\in H^1(\mathcal{D}_{\delta}):w|_{\ssphere}=0\right\rbrace.
$$
We seek a weak solution restricted to $H_{0, \delta}^1(\mathcal{D}_{\delta})$, for that we write the variational formulation as: Find $w\in H_{0, \delta}^1(\mathcal{D}_{\delta})$ such that:
\begin{alignat}{3}\label{eq:var_form}
\int_{\mathcal{D}_{\delta}}\tilde{\varepsilon}(\bx)\nabla w(\bx)\nabla \tilde{w}(\bx) &+ \int_{\mathcal{D}_{\delta}}\tilde{\lambda}(\bx)\widetilde{\mathcal{F}}\left(\overline{w}^{\hat{u}_1}(\bx)\right)w(\bx)\tilde{w}(\bx)\nonumber \\
&+\int_{\partial B_{\delta}(\bo)}\left(\mathcal{T}w\right)\tilde{w}(\bx) = \int_{\mathcal{D}_{\delta}}\tilde{f}(\bx)
\tilde{w}(\bx)\quad &&\forall\ \ \tilde{w}\in H_{0, \delta}^1(\mathcal{D}_{\delta}),
\end{alignat}
where we used $\tilde{\varepsilon}(\bx)=1$ on $\partial B_{\delta}(\bo)$ and $\mathcal{T}$ is the Dirichlet-to-Neumann operator of the harmonic extension in $B_{\delta}(\bo)$. 

Assume that we have an expansion of the Dirichlet boundary condition, then $w$ is restricted to $\partial B_{\delta}(\bo)$ as follows
\begin{equation}\label{eq:dirichlet_cond_1}
w|_{\partial B_{\delta}(\bo)}(\delta, \theta, \varphi) = \sum_{\ell, m}\gamma_{\ell m}Y_{\ell}^m(\theta, \varphi),\quad 0\leq\theta\leq \pi;\quad 0\leq \varphi\leq 2\pi.
\end{equation}
Then, we can extend $w$ harmonically by
\begin{equation}\label{eq:hamonic_extension}
w|_{B_{\delta}(\bo)}(r,\theta,\varphi) = \sum_{\ell, m}\gamma_{\ell m}\left(\frac{r}{\delta}\right)^{\ell}Y_{\ell}^m(\theta, \varphi),\quad 0\leq r\leq \delta;\quad0\leq\theta\leq \pi;\quad 0\leq \varphi\leq 2\pi.
\end{equation}
Let $\bn_{\delta}$ be the unit normal vector pointing outward on the ball $\partial B_{\delta}(\bo)$ with respect to the ball $B_{\delta}(\bo)$. As a consequence we compute the normal derivative $\partial_{\bn_{\delta}}w=\nabla w\cdot \bn_{\delta}$ on $\partial B_{\delta}(\bo)$ as
\begin{alignat}{3}\label{eq:neumann_harmonic}
\left( \mathcal{T}w\right)|_{B_{\delta}(\bo)}(\delta,\theta,\varphi) &=\partial_{\bn_{\delta}}w(\delta, \theta, \varphi)\nonumber \\
& =  \sum_{\ell, m}\gamma_{\ell m}\left(\frac{\ell}{\delta}\right)Y_{\ell}^m(\theta, \varphi),\qquad &&0\leq\theta\leq \pi;\quad 0\leq \varphi\leq 2\pi.
\end{alignat}

\begin{remark}
The bilinear form on the left side of Eq.~\eqref{eq:var_form},
\begin{equation}\label{eq:bilinear}
\begin{aligned}
a(w, \tilde{w}) & = \int_{\mathcal{D}_{\delta}}\varepsilon(\bx)\nabla w(\bx)\nabla \tilde{w}(\bx)+\int_{\mathcal{D}_{\delta}}\tilde{\lambda}(\bx) \widetilde{\mathcal{F}}\left(\overline{w}(\bx)\right) w(\bx)\tilde{w}(\bx) \\
& \quad+ \int_{\partial B_{\delta}(\bo)}\left( \mathcal{T} w(\bx)\right)\tilde{w}(\bx),
\end{aligned}
\end{equation}
is positive definite and symmetric by the definition of $\widetilde{\mathcal{F}}$ and the Dirichlet-to-Neumann operator $\mathcal{T}$ for given $w(\bx)$.
\end{remark}

\subsubsection{Galerkin Solution in Unit Ball}
For finding the functions belonging to $H_{0, \delta}^1(\mathcal{D}_{\delta})$; we introduce the radial functions
$$
\varrho_i(r) = (1-r)L_i'\left(\frac{2(r-\delta)}{1-\delta}-1\right),
$$
i.e., $\varrho_i(1) = 0$. Here $L_i$ denotes the Legendre polynomial of the $i^{\ith}$ degree. We then discretize both the radial part and the spherical part of the unknown $w$ by linear combination of the basis element $\lbrace\varrho_iY_{\ell}^m\rbrace$ with $1\leq i\leq N;0\leq \ell\leq \ell_{\max};$ and $-\ell\leq m\leq \ell$, where $N$ denotes the maximum degree of Legendre polynomials and $\ell_{\max}$ denotes the maximum number of spherical harmonics. We denote the space spanned by these elements by $\mathcal{B}_{N,\ell_{\max}}(\mathcal{D}_{\delta})$ which is defined as
\begin{alignat*}{3}
\mathcal{B}_{N,\ell_{\max}}(\mathcal{D}_{\delta})&=\mathrm{span}\left\lbrace \varrho_i(r)Y_{\ell}^m(\theta, \varphi):1\leq i\leq N,\quad 0\leq\ell\leq  \ell_{\max},\quad -\ell\leq m\leq \ell\right\rbrace\nonumber \\
&\subset H_{0, \delta}^1(\mathcal{D}_{\delta}).
\end{alignat*}
Then, the Galerkin discretisation of the variational formulation~\eqref{eq:var_form} reads: Find $w_{\mathcal{B}}\in\mathcal{B}_{N,\ell_{\max}}(\mathcal{D)}$, such that
\begin{equation}\label{eq:galerkin_temp_1}
a\left(w_{\mathcal{B}},\vry\right)=\int_{\mathcal{D}_{\delta}}\tilde{f}\vry,\quad\forall\ \ \vry\in\mathcal{B}_{N,\ell_{\max}}(\mathcal{D}_{\delta}),
\end{equation}
where $a(\cdot,\cdot)$ is given by Eq.~\eqref{eq:bilinear}. Since $w_{\mathcal{B}}\in \mathcal{B}_{N,\ell_{\max}}(\mathcal{D}_{\delta})$, we can write $w_{\mathcal{B}}$ in the form
\begin{equation}\label{eq:galerkin_temp_2}
w_{\mathcal{B}}(r,\theta, \varphi) = \sum_{i=1}^N\sum_{\ell, m}\left[\phi_{\mr}\right]_{i\ell}^m\varrho_i(r)Y_{\ell}^m(\theta,\varphi),\quad \forall\ \ \delta\leq r\leq 1;\ 
 0\leq \theta\leq \pi;\ \ 0\leq \varphi \leq 2\pi,
\end{equation}
and consequently
\begin{equation}\label{eq:galerkin_temp_3}
\left(\mathcal{T}w_{\mathcal{B}}\right)|_{B_{\delta}(\bo)}(\delta ,\theta, \varphi) = \sum_{i=1}^N\sum_{\ell, m}\left[\phi_{\mr}\right]_{i\ell}^m\left( \frac{\ell}{\delta}\right)\varrho_i(\delta)Y_{\ell}^m(\theta,\varphi),\ \forall\ 0\leq \theta\leq \pi;\ 0\leq \varphi \leq 2\pi,
\end{equation}
where $\left[\phi_{\mr}\right]_{i\ell}^m$ is the real coefficient of $w_{\mathcal{B}}$ corresponding to the node $\varrho_iY_{\ell}^m$. 

Substituting Eq.~\eqref{eq:galerkin_temp_2} and Eq.~\eqref{eq:galerkin_temp_3} in Eq.~\eqref{eq:galerkin_temp_1} and taking the test function $\vry=\varrho_j(r)Y_{\ell'}^{m'}(\theta, \varphi)$, we then obtain a system of equations for all $1\leq j \leq N,0\leq \ell'\leq \ell_{\max},$ and $-\ell'\leq m'\leq \ell'$
\begin{alignat}{3}\label{eq:system_equation_1}
\sum_{i=1}^N\sum_{\ell, m}\left[\phi_{\mr}\right]_{i\ell}^m\Bigg(&\int_{\mathcal{D}_{\delta}}\tilde{\varepsilon}(\bx)\nabla \Big( \varrho_iY_{\ell}^m\Big)\cdot \nabla\Big(\varrho_jY_{\ell'}^{m'}\Big) + \int_{\mathcal{D}_{\delta}}\tilde{\lambda}(\bx) \widetilde{\mathcal{F}}\left(\overline{w}_{\mathcal{B}}^{\hat{u}_1}(\bx)\right)\varrho_iY_{\ell}^m\varrho_jY_{\ell'}^{m'} \nonumber \\
&+\frac{\ell}{\delta}\int_{\partial B_{\delta}(\bo)}\varrho_iY_{\ell}^m\varrho_jY_{\ell'}^{m'}\Bigg) =  \int_{\mathcal{D}_{\delta}}\tilde{f}\varrho_jY_{\ell'}^{m'}.
\end{alignat}
In order to write a system of equations, we define the index
$$
k=N\left( \ell^2+m+1\right)+i\in \lbrace 1,2,\dots,N\left(\ell_{\max}+1\right)^2\rbrace,
$$
which corresponds to the triple $(i,\ell,m)$. Let $k$ corresponds to $(i,\ell,m)$ and $k'$ corresponds to $(j,\ell',m')$. Then Eq.~\eqref{eq:system_equation_1} can be recast as
\begin{equation}\label{eq:system_equation_2}
\overline{\bA}\ \overline{X}_{\mr}=\overline{\bF},
\end{equation}
where $\overline{\bA}$ is a matrix of dimension $N(\ell_{\max}+1)^2\times N(\ell_{\max}+1)^2$ with elements $\left(\overline{\bA}\right)_{k,k'}$ for all $1\leq k,k'\leq N(\ell_{\max}+1)^2$, defined by
\begin{alignat}{3}\label{eq:matrix_A_entry}
\left(\overline{\bA}\right)_{k,k'} =& \int_{\mathcal{D}_{\delta}}\tilde{\varepsilon}(\bx)\nabla \left( \varrho_iY_{\ell}^m\right)\cdot \nabla\left(\varrho_jY_{\ell'}^{m'}\right)+ \int_{\mathcal{D}_{\delta}}\tilde{\lambda}(\bx) \widetilde{\mathcal{F}}\left(\overline{w}_{\mathcal{B}}^{\hat{u}_1}(\bx)\right)\varrho_iY_{\ell}^m\varrho_jY_{\ell'}^{m'}\nonumber \\
& + \frac{\ell}{\delta}\int_{\partial B_{\delta}(\bo)}\varrho_iY_{\ell}^m\varrho_jY_{\ell'}^{m'},
\end{alignat}
$\overline{X}_{\mr}$ is the column vector of $N\left( \ell_{\max}+1\right)^2$ unknowns $ \left[\phi_{\mr}\right]_{i\ell}^m$, i.e.,
\begin{equation}\label{eq:vector_x_entry}
(\overline{X}_{\mr})_{k} = \left[\phi_{\mr}\right]_{i\ell}^m,\quad \forall\ \ k\in \lbrace 1,\dots, N(\ell_{\max}+1)^2\rbrace,
\end{equation}
and $\overline{\bF}$ is a column vector with $N(\ell_{\max}+1)^2$ entries defined by
\begin{equation}\label{eq:vector_F_entry}
(\overline{\bF})_{k'} = \int_{\mathcal{D}_{\delta}}\tilde{f}\varrho_jY_{\ell'}^{m'},\quad \forall\ \ k'\in \lbrace 1,\dots, N(\ell_{\max}+1)^2\rbrace.
\end{equation}

To summarize, in order to solve Eq.~\eqref{eq:var_form} we need to solve Eq.~\eqref{eq:system_equation_2} to obtain $\left[\phi_{\mr}\right]_{i\ell}^m$ and then obtain an approximate solution $w_{\mathcal{B}}(r,\theta, \varphi)\in \mathcal{B}_{N,\ell_{\max}}(\mathcal{D}_{\delta})$ according to Eq.~\eqref{eq:galerkin_temp_2}. Since, $w_{\mathcal{B}}$ is harmonic in $B_{\delta}(\bo)$, $w_{\mathcal{B}}$ can be extended harmonically in the ball $B_{\delta}(\bo)$ following Eq.~\eqref{eq:hamonic_extension} and hence we obtain an approximate solution defined in $B_1(\bo)$ to Eq.~\eqref{eq:harmonic_eq_2}.

\begin{remark}\label{rem:intergtation_torus}
The final thing remaining in the discretization process is the evaluation of the integrals in $\overline{\bA}$ and $\overline{\bF}$. We have integrals over the annulus $\mathcal{D}_{\delta}$ and $\partial B_{\delta}(\bo)$. Here, we follow the ideas from \cite{QSM18}; for completeness, we present a brief overview. The integral over $\partial B_{\delta}(\bo)$ is given by
\begin{equation*}
\begin{aligned}
\frac{\ell}{\delta}\int_{\partial B_{\delta}(\bo)}\varrho_i Y_{\ell}^m\varrho_j Y_{\ell'}^{m'} &= \ell\delta \varrho_i(\delta)\varrho_j(\delta)\int_{\mathbb{S}^2}Y_{\ell}^m Y_{\ell'}^{m'}\\
& = \ell \delta \varrho_i(\delta)\varrho_j(\delta)\delta_{\ell \ell'}\delta_{m m'},
\end{aligned}
\end{equation*}
where $\delta_{\ell \ell'}$ and $\delta_{m m'}$ are the Kronecker delta functions.

The integral over $\mathcal{D}_{\delta}$ can be divided into radial and spherical parts. Let $h\in L^1\left(B_1(\bo)\right)$ then the integral of $h$ over $\mathcal{D}_{\delta}$ can be written separately as
$$
\int_{\mathcal{D}_{\delta}}h(\bx)d\bx=\int_{\delta}^1r^2\int_{\mathbb{S}^2}h(r,\bs)d\bs dr,\quad \bs\in \mathbb{S}^2,
$$
and $\bx=r\bs$. The spherical part can be computed using the Lebedev quadrature \cite{Hax07}. For the radial part we use the Legendre-Gau\ss-Lobatto (LGL) quadrature rule \cite{Par99} defined by quadrature points $x_{m}\in [-1,1]$ and the quadrature weights $\omega_m^{\lgl}$, $1\leq m\leq N_{\lgl}$ for $N_{\lgl}$ quadrature points. Using the change of variable
$$
r=\tfrac{1-\delta}{2}\left(x+1\right)+\delta,\quad x\in [-1,1],
$$
we approximate the integral by the following quadrature rule
\begin{alignat*}{3}
\int_{\mathcal{D}_{\delta}}h(\bx)d\bx \approx & \tfrac{1-\delta}{2}\sum_{m=1}^{N_{\lgl}}\sum_{n=1}^{N_{\leb}}\omega_m^{\lgl}\omega_n^{\leb}\left(\tfrac{1-\delta}{2}(x_m+1)+\delta\right)^2 h\left(\tfrac{1-\delta}{2}(x_m+1)+\delta,\bs_n\right).
\end{alignat*}
\end{remark}

\section{Numerical Simulations}\label{sec:numres}
This section presents some numerical studies for the ddNPB-SES method and compares them with the ddLPB-SES method.

We first introduce the electrostatic solvation energy for the NPB equation. For the NPB, the solvation energy becomes more involved than the LPB. We follow the definition of the energy from \cite{SH90}, which is given by
\begin{alignat*}{2}
E_s&=\frac{\beta}{2}\int_{\Omega_0}\rho^{\sol}(\bx)\psi_{\mr}(\bx)\nonumber \\
&\ \ + \frac{\beta^2\kappa^2\varepsilon_s}{8\pi} \int_{\Omega_0}\lambda(\bx)\Big[\psi_{\mr}(\bx)\sinh\left(\psi_{\mr}(\bx)\right) -2\left\lbrace\cosh\left(\psi_{\mr}(\bx)\right)-1\right\rbrace\Big].
\end{alignat*}

\begin{remark}
The second term in $E_s$ is referred to as the electrostatic stress, and the last term is the osmotic pressure, see \cite{SH90}. In the case of the LPB equation, we can linearize the last two terms and get $E_s$ corresponding to the ddLPB method (see \cite{QSM19}).
\end{remark}

Now, we provide more details about solving the system of nonlinear equations. We notice that in matrix $\overline{\bA}$, the first and third parts in Eq.\eqref{eq:matrix_A_entry} are constant throughout the iterative process, and hence they need to be computed only once and can be used later. For computing the solution at the $(\nu)^{\ith}$ iterative step, we follow a damping approach in which
$$
X^{(\nu)} = X^{(\nu-1)}+\omega\left(X^{\mathrm{aux}}-X^{(\nu-1)}\right),\qquad \nu\geq 1
$$
where $X^{\mathrm{aux}}$ is obtained by solving Eq.~\eqref{eq:system_equation_2} and $\omega\in (0,1]$ is a damping parameter. One can choose a constant or an adaptive damping parameter. For nonlinear iterations, an adaptive choice of $\omega$ yields better results (see \cite{SG90}), but in this work, the constant damping gave satisfactory results. Hence, we fix $\omega=0.5$ for the simulations, which is obtained empirically.

We have three different iteration loops in our method. We refer to the global iteration process as outer iterations (indexed by $\tk$), and the convergence is reached when the relative increment satisfies
$$
\mathtt{inc}_{\tk}:=\frac{|E_s^{\tk}-E_s^{\tk-1}|}{|E_s^{\tk}|}\leq \mathtt{tol}, \qquad \mathrm{for}\ \ \tk\geq 1,
$$
and given tolerance, \texttt{tol}. The second is the dd-iterative loop for the GSP solver (Eq.~\eqref{eq:ndr_subdomain}--\eqref{eq:ndr_bc}); here, the stopping criteria is the relative reduction of the reaction potential but with a relaxed stopping criteria of $10\times \mathtt{tol}$. Lastly is the nonlinear loop for solving Eq.~\eqref{eq:system_equation_2}. Here we again use the relative reduction of the solution vector with the stopping criteria of $100\times \mathtt{tol}$. Along with this stopping criteria, we also have a maximal number of iterative loops for each iteration. To solve the outer iterative loop, we start with the zero solution as the initial iterate, but we use the solution from the previous iterative step for subsequent iterations. We follow the same procedure for the other iterative loops as well. Unless specified by default, we set $\mathtt{tol}=10^{-6}$ and solve the system of equations using the LU decomposition method.

Lastly, we would like to mention the constants arising in the method. By default, we assume the solute cavity in a vacuum and the solvent to be water. Hence, the relative dielectric permittivity of the solute cavity is one and $\varepsilon_s=78.54$ at room temperature, $T=298.15$ K. Further, we set the Debye H\"uckel constant, $\kappa=0.104$\ \AA$^{-1}$ for an ionic strength of $0.1$ molar. We use the Hartree energy unit system and thus $4\pi \varepsilon_\abso=1$. We read the input files in\ \AA\ units, but then we internally convert them to the Hartree units, and hence, the distance is represented by atomic units (a.u.). The atomic centers, charges,  and the vdW radii are obtained from the \textsc{PDB} files \cite{Berman00} and the \textsc{PDB2PQR} package \cite{DCLN07, Dolinsky04}. Finally, we set the probe radius $r_p=1.4$\ \AA. All the simulations were performed on the in-house \textsc{MATLAB}.

\subsection{GSP Solver}\label{subsec:gsp_solver}
We first present results for the GSP solver presented in Sec.~\ref{sec:gsp_solver} in a unit ball. Unless mentioned, we assume a $0.1$ charge at the origin with a vdW radii of $2$, i.e., $r_1=2$\ \AA.  We also assume the absence of the Stern layer and hence set $a=0$\ \AA. As we have a single atom, we have the rotational symmetry of the system. Hence, we define the dielectric permittivity and the ion-exclusion function only as a radial variable,  i.e.,
\begin{alignat*}{3}
\varepsilon(r)&=\left\lbrace\begin{array}{ll}
1& r<r_1,\\
1+(\varepsilon_s-1)\xi\left( \dfrac{r-r_1}{r_p}\right)& r_1\leq r\leq r_1+r_p, \\
\varepsilon_s & r>r_1+r_p.
\end{array}\right.,\nonumber \\
\lambda(r)&=\left\lbrace \begin{array}{ll}
0& r<r_1,\\
\xi\left( \dfrac{r-r_1}{r_p}\right)& r_1\leq r\leq r_1+r_p,\\
1 & r>r_1+r_p.
\end{array}\right..
\end{alignat*}

First, we plot various potentials with respect to the radial component for this system, Fig.~\ref{fig:electrostatic_potential_single_atom}. Here, we set $r_0=1$\ \AA\ and the radial discretization parameters to $N=20$ and $N_{\lgl}=200$. The Bessel extension is obtained by extending $\psi_{\me}$ on the outer boundary using the modified spherical Bessel function of the second kind \cite[Chapter~14]{AWH13}, and the harmonic extension is obtained by extending $w_{\mathcal{B}}$ (see Eq.~\eqref{eq:galerkin_temp_2}) on the inner boundary. We notice that we have both extensions' continuity and a zero jump of the normal derivative as predicted by the theory.

\input{figures/potential_one_sphere.tex}

Next, we show the variation between the potentials obtained by using the LPB and the NPB equation. For this, we define a function, $\mathrm{Var}_\psi$, given by
$$
\mathrm{Var}_\psi(r):=\left|\psi_{\mathrm{NPB}}(r)-\psi_{\mathrm{LPB}}(r)\right|\qquad \forall r_1\leq r\leq R_1,
$$
where $\psi_{\mathrm{NPB}}$ and $\psi_{\mathrm{LPB}}$ refer to the reaction potential of the NPB and the LPB equation, respectively.  As, $r_1$ and $r_p$ are fixed, the variation in $r$ comes from $r_0$. We expect that after a certain $r_0$,  $\mathrm{Var}_\psi$ should decrease. For this simulation, we set $r_0=10,000$\ \AA\ and the radial discretization parameters $N=50$ and $N_{\lgl}=5000$.  We take three different charges, $q_1\in \lbrace 10^{-1}, 10^{-2}, 10^{-4}\rbrace$. In Fig.~\ref{fig:plot_variation}, we notice that for $q_1=10^{-1}$ after $r>10^3$ the $\mathrm{Var}_\psi\leq 10^{-3}$. As the atom's charge decreases, the same tolerance can be achieved for a smaller $r$. This example shows the importance of having a nonlinear region in the NPB equation, especially for a highly charged molecule, as the potential obtained by the NPB equation does not satisfy the low-potential condition near the SES and thus the LPB is not accurate in this region.

\input{figures/var_psi.tex}

\subsection{Convergence of Global Strategy}\label{subsec:convergence_global}
In this example, we consider the caffeine molecule to show the convergence of the global strategy introduced in Sec.~\ref{sec:global_strategy}. Caffeine having twenty-four atoms gives a good idea about the developed method.

The Schwarz domain decomposition method, which is used to solve Eq.~\eqref{eq:ndr_subdomain} and Eq.~\eqref{eq:hsp_subdomain}, is well studied and its convergence can be guaranteed \cite{RS21} in a continuous setting. To study the convergence of our global strategy, we set the discretization parameters as $\ell_{\max}=9$, $N_\leb=350$, $N=15$, and $N_\lgl=30$. With this, we compute an ``exact" solvation energy, $E_s^\infty$ for 15 outer iterations and then define an error function as
$$
\mathtt{Error}_{N_{\mathrm{it}}}:=|E_s^\infty -E_s^{N_{\mathrm{it}}}|,
$$
where $N_{\mathrm{it}}$ are the number of outer iterations.  The geometric parameters are fixed as $r_0=5$\ \AA, and $a=1$\ \AA.

\begin{figure}[!t]
\centering
\begin{tikzpicture}[scale=0.5]
\begin{axis}[
    legend pos=north east, xlabel = \Large{$N_{\mathrm{it}}$}, ylabel= \Large{$E_s$ (Hartree)},     legend cell align ={left}, title = \large{Solvation Energy}]
\addplot[color=dark_green,  mark=oplus*, line width = 0.5mm, dashdotted,,mark options = {scale= 1.5, solid}]
coordinates{( 1 , 0.010421279173 )
( 2 , 0.008593437004 )
( 3 , 0.008485205399 )
( 4 , 0.008584746765 )
( 5 , 0.008554823368 )
( 6 , 0.00852163028 )
( 7 , 0.008573949897 )
( 8 , 0.008523791078 )
( 9 , 0.008562530425 )
( 10 , 0.008535403728 )
( 11 , 0.008552681673 )
( 12 , 0.008542407509 )
( 13 , 0.008548157283 )
( 14 , 0.008545107949 )
( 15 , 0.008546612261 )
};
\end{axis}
\end{tikzpicture}\hspace*{0.1em}
\begin{tikzpicture}[scale=0.5]
\begin{semilogyaxis}[
    legend pos=north east, xlabel = \Large{$N_{\mathrm{it}}$}, ylabel= \Large{$\mathtt{Error}_{N_{\mathrm{it}}}$},     legend cell align ={left}, title = \large{\texttt{Error}}]
\addplot[color=copper,  mark=oplus*, line width = 0.5mm, dashdotted,,mark options = {scale= 1.5, solid}]
coordinates{( 1 , 0.0018746669120000004 )
( 2 , 4.682474300000093e-05 )
( 3 , 6.140686199999912e-05 )
( 4 , 3.813450400000015e-05 )
( 5 , 8.211107000000162e-06 )
( 6 , 2.4981981000000292e-05 )
( 7 , 2.7337635999999735e-05 )
( 8 , 2.282118299999983e-05 )
( 9 , 1.5918163999999527e-05 )
( 10 , 1.1208533000000437e-05 )
( 11 , 6.0694120000009955e-06 )
( 12 , 4.204751999999437e-06 )
( 13 , 1.545022000001145e-06 )
( 14 , 1.5043119999998827e-06 )
( 15 , 0.0 )
};
\end{semilogyaxis}
\end{tikzpicture}\\
\begin{tikzpicture}[scale=0.5]
\begin{axis}[
    legend pos=north east, xlabel = \Large{$N_{\mathrm{it}}$}, ylabel= \Large{$N_{\mathrm{dd}}$},     legend cell align ={left}, title = \large{Number of dd-loops}]
\addplot[color=carrotorange,  mark=oplus*, line width = 0.5mm, dashdotted,,mark options = {scale= 1.5, solid}]
coordinates{( 1 , 15.0 )
( 2 , 15.0 )
( 3 , 15.0 )
( 4 , 15.0 )
( 5 , 15.0 )
( 6 , 13.0 )
( 7 , 11.0 )
( 8 , 11.0 )
( 9 , 11.0 )
( 10 , 11.0 )
( 11 , 10.0 )
( 12 , 9.0 )
( 13 , 8.0 )
( 14 , 7.0 )
( 15 , 6.0 )
};
\end{axis}
\end{tikzpicture}
\caption{Example~\ref{subsec:convergence_global}: Electrostatic solvation energy (top left), error (top right), and number of dd loops (bottom) for the caffeine molecule with respect to the number of outer iterations. 
}\label{fig:error_global_convergece}
\end{figure}
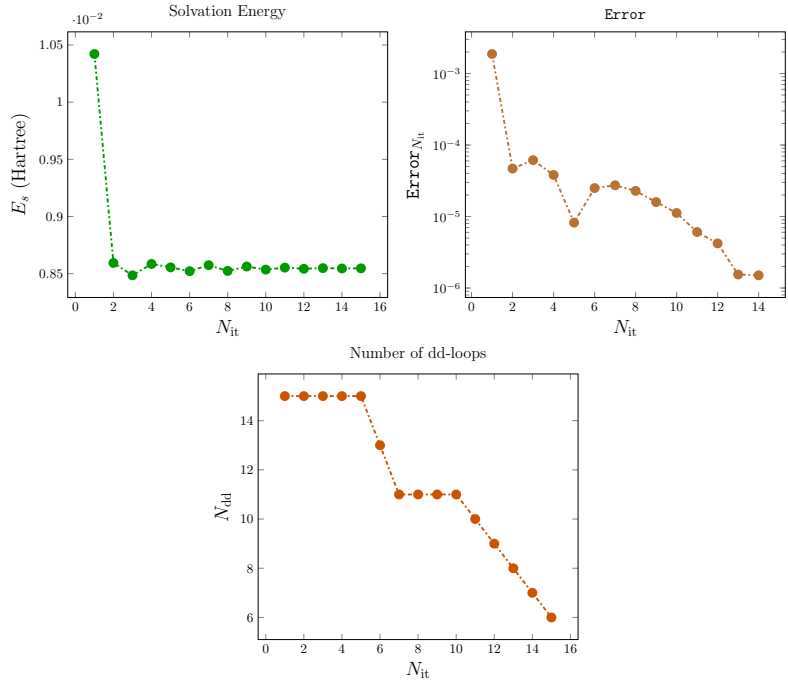

We observe in Fig.~\ref{fig:error_global_convergece} that the $E_s$ converges for $N_{\mathrm{it}}$ (top left), and the process stops when the desired tolerance is reached. The error also decreases monotonically with increasing $N_{\mathrm{it}}$ (top right). We finally present the number of domain decomposition ($N_{\mathrm{dd}}$) loops required to solve Eq.~\eqref{eq:ndr_subdomain} in the global iterative process (bottom). As the number of outer loops increases, the number of dd loops monotonically decreases. This example also gives a good idea for choosing the maximum outer iterations, around $5$ in this case.

\subsection{Effect of Discretisation Parameters}\label{subsec:disc_param}
In this example, we study the effects of the discretization parameters for the formaldehyde molecule.

The geometric parameters are set as $r_0=3$\ \AA\ and $a=1$\ \AA. For the reference (``exact'') solvation energy we use the discretisation parameters $\ell_{\max}=11$, $N_{\leb}=1202$, $N=15$, and $N_{\lgl}=50$. 

In Fig.~\ref{fig:formal_error_spherical} the $\ell_{\max}$ is varied from $2$ to $11$ (left)  and $N_{\leb}$ is varied from $194$ to $1202$ (right), while $N=15$ and $N_{\lgl}=50$. We can observe that the proposed algorithm improves systematically with an increase in the number of parameters. The results for the radial discretization are presented in Fig.~\ref{fig:formal_error_radial} where $N$ is varied from $5$ to $15$ (left) and $N_{\lgl}$ is varied from $15$ to $50$ (right), while $\ell_{\max}=11$ and $N_{\leb}=1202$. We have similar observations as compared to the spherical discretization parameters, i.e., the algorithm improves systematically. This example gives a good choice in selecting the discretization parameters.

\begin{figure}[!t]
\centering
\begin{tikzpicture}[scale=0.5]
\begin{axis}[
    legend pos=north east, xlabel = \Large{$\ell_{\max}$}, ylabel= \Large{$E_s$\ (Hartree)},  ymin= 0.000117262430, ymax= 0.00033248063,     legend cell align ={left}, title = \large{Solvation Energy}]
\addplot[color=crimson,  mark=oplus*, line width = 0.5mm, dashdotted,, mark options = {scale= 1.5, solid}] 
coordinates{
(3.000000000000, 0.000137262430)
(4.000000000000, 0.000227953777)
(5.000000000000, 0.000223316824)
(6.000000000000, 0.000216332912)
(7.000000000000, 0.000211933053)
(8.000000000000, 0.000214700060)
(9.000000000000, 0.000216097272)
(10.000000000000, 0.000199927560)
(11.000000000000, 0.000224871539)
};
\addlegendentry{\large{$N_{\leb}=1202$}} 
\addplot[color=black,  line width = 0.5mm] 
coordinates{( 2,  0.000224871539)( 12, 0.000224871539)};
\addlegendentry{\large{$\ell_{\max}=11, N_{\leb}=1202$}} 
\end{axis}
\end{tikzpicture}\hspace{1em}
\begin{tikzpicture}[scale=0.5]
\begin{axis}[
    legend pos=north east, xlabel = \Large{$N_{\leb}$}, ylabel=  \Large{$E_s$\ (Hartree)}, ymin= 0.00010147515, ymax= 0.000348267925,  legend cell align ={left}, title = \large{Solvation Energy}]
\addplot[color=awesome,  mark=oplus*, line width = 0.5mm, dashdotted,, mark options = {scale= 1.5, solid}] 
coordinates{
(194.000000000000, 0.000194736057)
(230.000000000000, 0.000328267925)
(266.000000000000, 0.000210018865)
(302.000000000000, 0.000227937517)
(350.000000000000, 0.000200268527)
(434.000000000000, 0.000195294730)
(590.000000000000, 0.000225944791)
(770.000000000000, 0.000211751187)
(974.000000000000, 0.000218439971)
(1202.000000000000, 0.000224871539)
};
\addlegendentry{\large{$\ell_{\max}=11$}} 
\addplot[color=black,  line width = 0.5mm] 
coordinates{( 100 ,  0.000224871539)( 1350,  0.000224871539)};
\addlegendentry{\large{$\ell_{\max}=11, N_{\leb}=1202$}} 
\end{axis}
\end{tikzpicture}
\caption{Example~\ref{subsec:disc_param}: Electrostatic solvation energy contribution with respect to $\ell_{\max}$ (left) by setting $N_{\leb}=1202$, and with respect to $N_{\leb}$ (right) by setting $\ell_{\max}=11$ for the formaldehyde molecule. The ``exact'' solvation energy (solid line) is obtained with $\ell_{\max}=11$, $N_{\leb}=1202$, $N=15$, and $N_{\lgl}=50$.
}\label{fig:formal_error_spherical}
\end{figure}
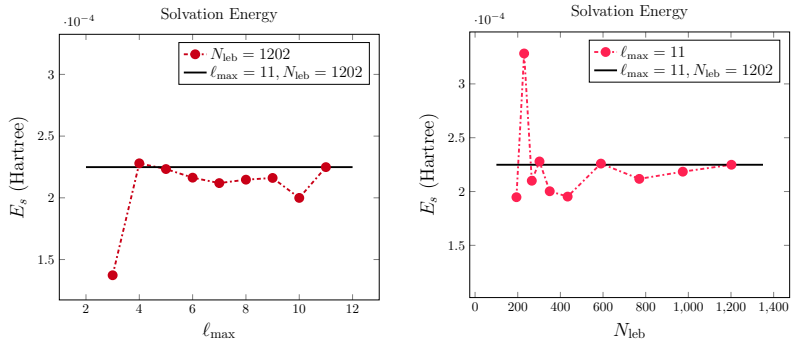

\begin{figure}[!t]
\centering
\begin{tikzpicture}[scale=0.5]
\begin{axis}[
    legend pos=north east, xlabel = \Large{$N_{\textcolor{white}{aj}}$}, ylabel= \Large{$E_s$\ (Hartree)},  ymin= -0.00002267118, ymax= 0.00049041425,     legend cell align ={left}, title = \large{Solvation Energy}]
\addplot[color=blue,  mark=oplus*, line width = 0.5mm, dashdotted,, mark options = {scale= 1.5, solid}] 
coordinates{
(5.000000000000, 0.000119328814)
(6.000000000000, 0.000203365118)
(7.000000000000, 0.000239246490)
(8.000000000000, 0.000183627761)
(9.000000000000, 0.000172605395)
(10.000000000000, 0.000228958354)
(11.000000000000, 0.000192673221)
(12.000000000000, 0.000228577573)
(13.000000000000, 0.000232312276)
(14.000000000000, 0.000200702721)
(15.000000000000, 0.000224871539)
};
\addlegendentry{\large{$N_{\lgl}=50$}} 
\addplot[color=black,  line width = 0.5mm] 
coordinates{( 3,  0.000224871539 )( 17, 0.000224871539)};
\addlegendentry{\large{$N=15, N_{\lgl}=50$}} 
\end{axis}
\end{tikzpicture}\hspace{1em}
\begin{tikzpicture}[scale=0.5]
\begin{axis}[
    legend pos=north east, xlabel = \Large{$N_{\lgl}$}, ylabel=  \Large{$E_s$\ (Hartree)}, ymin=0.000025732370, ymax= 0.00042401069,  legend cell align ={left}, title = \large{Solvation Energy}]
\addplot[color=royal_blue,  mark=oplus*, line width = 0.5mm, dashdotted,, mark options = {scale= 1.5, solid}] 
coordinates{
(15.000000000000, 0.000185732370)
(20.000000000000, 0.000217993632)
(25.000000000000, 0.000202618389)
(30.000000000000, 0.000195031796)
(35.000000000000, 0.000229166497)
(40.000000000000, 0.000217183857)
(45.000000000000, 0.000194482307)
(50.000000000000, 0.000224871539)
};
\addlegendentry{\large{$N=15$}} 
\addplot[color=black,  line width = 0.5mm] 
coordinates{( 10, 0.000224871539)(55,  0.000224871539)};
\addlegendentry{\large{$N=15, N_{\lgl}=50$}} 
\end{axis}
\end{tikzpicture}
\caption{Example~\ref{subsec:disc_param}: Electrostatic solvation energy contribution with respect to $N$ (left) by setting $N_{\lgl}=50$, and with respect to $N_{\lgl}$ (right) by setting $N=15$ for the formaldehyde molecule. The ``exact'' solvation energy (solid line) is obtained with $\ell_{\max}=11$, $N_{\leb}=1202$, $N=15$, and $N_{\lgl}=50$.
}\label{fig:formal_error_radial}
\end{figure}

\begin{figure}[t!]
\centering
\includegraphics[scale=0.15]{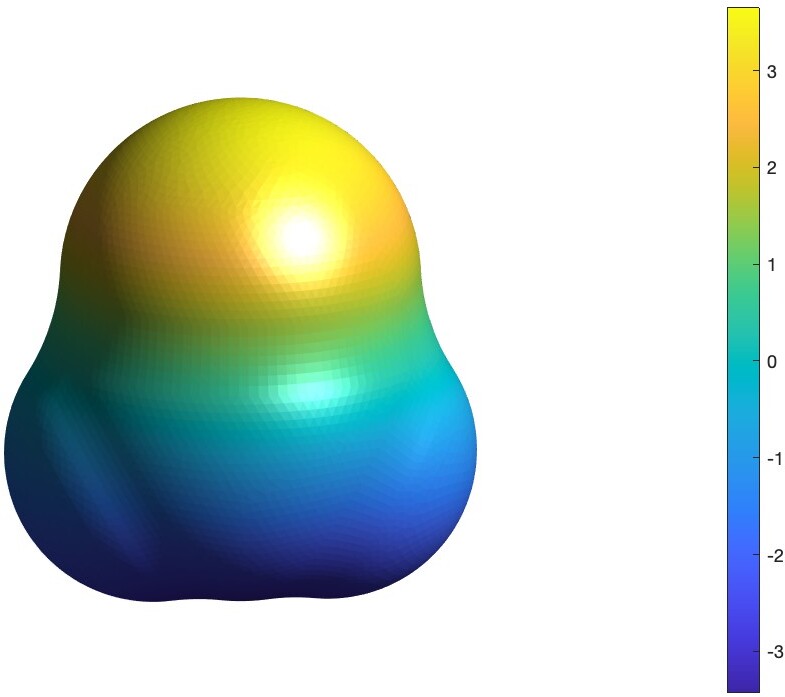}
\caption{Example~\ref{subsec:disc_param}: Reaction potential for the formaldehyde molecule on $\Omega_{\ses}$. 
Discretisation parameters set to $\ell_{\max}=11$, $N_{\leb}=1202$, $N=15$, and $N_{\lgl}=50$. 
}\label{fig:reaction_potential_ch4}
\end{figure}

Fig.~\ref{fig:reaction_potential_ch4} shows the electrostatic potential on the SES surface for this example. We observe the mirror symmetry for this molecule.

\subsection{Stern Layer Length}\label{subsec:stern_layer}
Until now, we have not discussed the effects of the Stern layer length. In this example, we study its effect. We consider the hydrogen fluoride molecule with $r_0=2$\ \AA. In this example we set the discretization parameters as $\ell_{\max}=9$, $N_{\leb}=350$, $N=15$, and $N_{\lgl}=50$.

\begin{figure}[!t]
\centering
\begin{tikzpicture}[scale=0.5, spy using outlines={rectangle, magnification=5,connect spies}]
\begin{semilogyaxis}[
    legend pos=north east, xlabel = \Large{$a$\ (a.u.)}, ylabel= \Large{$\mathrm{Diff}_{E_s}$},  legend cell align ={left}]
\addplot[color=crimson,  mark=oplus*, line width = 0.5mm, dashdotted,,mark options = {scale= 1.5, solid}]
coordinates{(0.000000000000, 0.306691180603)
(1.000000000000, 0.001871523045)
(2.000000000000, 0.000237048190)
(3.000000000000, 0.000063300281)
(4.000000000000, 0.000024269297)
(5.000000000000, 0.000011509595)
(6.000000000000, 0.000006448935)
(7.000000000000, 0.000004170961)
(8.000000000000, 0.000002771662)
(9.000000000000, 0.000001966313)
(10.000000000000, 0.000001414599)
(11.000000000000, 0.000001063338)
(12.000000000000, 0.000000818138)
(13.000000000000, 0.000000632638)
(14.000000000000, 0.000000505152)
(15.000000000000, 0.000000409358)
(16.000000000000, 0.000000336446)
(17.000000000000, 0.000000280282)
(18.000000000000, 0.000000230331)
(19.000000000000, 0.000000194603)
(20.000000000000, 0.000000165381)
(21.000000000000, 0.000000140796)
(22.000000000000, 0.000000120566)
(23.000000000000, 0.000000103828)
(24.000000000000, 0.000000088819)
(25.000000000000, 0.000000077197)
(26.000000000000, 0.000000067027)
(27.000000000000, 0.000000058554)
(28.000000000000, 0.000000051371)
(29.000000000000, 0.000000044835)
(30.000000000000, 0.000000039358)
};
\end{semilogyaxis}
\end{tikzpicture}
\caption{Example~\ref{subsec:stern_layer}: Difference of the electrostatic stress and the osmotic pressure with respect to Stern layer, $a$ for the hydrogen fluoride molecule. 
}\label{fig:stern_layer}
\end{figure}
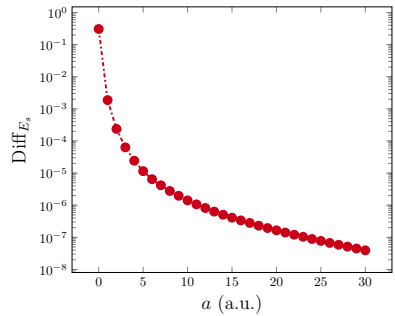

We vary $a$ from $0$ to $30$\ \AA. For $a=0$\ \AA, we have an absence of a Stern layer; hence, the ions are close to the $\ses$ surface. For this example we define a function $\mathrm{Diff}_{E_s}$ where
$$
\mathrm{Diff}_{E_s}=\Bigg|\frac{\beta^2\kappa^2\varepsilon_s}{8\pi} \int_{\Omega_0}\lambda(\bx)\Big[\psi_{\mr}(\bx)\sinh\left(\psi_{\mr}(\bx)\right) -2\left\lbrace\cosh\left(\psi_{\mr}(\bx)\right)-1\right\rbrace\Big]\Bigg|,
$$
i.e., the difference between the electrostatic stress and the osmotic pressure. In Fig.~\ref{fig:stern_layer} we plot $\mathrm{Diff}_{E_s}$ and, we notice that as $a$ increases, the $\mathrm{Diff}_{E_s}$ decreases and after $a>12$, the value is less than $10^{-6}$. One expects this observation as including the Stern layer leads to a screening effect near the charged surface. As a result, $\mathrm{Diff}_{E_s}$ is lower.

\subsection{Rotational Symmetry}\label{subsec:rot_sym}
One of the advantages of the dd-methods is the rotational invariance in the computation of $E_s$, \cite{QSM18, Jha23}. In this example, we show that the ddNPB-SES also satisfies this property. We consider the hydrogen fluoride molecule and we fix the hydrogen atom at the center and rotate the fluorine atom around the hydrogen atom in the $y-z$ plane.

The geometric parameters are set to $r_0=3$\ \AA\ and $a=5$\ \AA. The radial discretization parameters are set to $N=15$ and $N_{\lgl}=50$. We present results with respect to two sets of spherical discretization parameters, $\ell_{\max}=7$ and $11$. The number of Lebedev quadrature points is set according to \cite{CMS13} to get the exact quadrature.  For $\ell_{\max}=7$ we set $N_{\leb}=86$ and for $\ell_{\max}=11$, $N_{\leb}=1202$.

\begin{figure}[!t]
\centering
\begin{tikzpicture}[scale=0.5]
\begin{axis}[
    legend pos=north east, xlabel = \Large{$\theta$}, ylabel= \Large{$E_s$\ (Hartree)},  ymin= -0.000071396396, ymax= -0.000046953573,   legend cell align ={left}, title = \large{Solvation Energy}]
\addplot[color=crimson,  mark=oplus*, line width = 0.5mm, dashdotted,, mark options = {scale= 1.5, solid}] 
coordinates{
(0.000000000000, -0.000067056107)
(0.314200000000, -0.000053191823)
(0.628400000000, -0.000060808619)
(0.942600000000, -0.000058772433)
(1.256800000000, -0.000050953573)
(1.571000000000, -0.000065424914)
(1.885200000000, -0.000053145113)
(2.199400000000, -0.000059188722)
(2.513600000000, -0.000057872574)
(2.827800000000, -0.000052841925)
(3.142000000000, -0.000067056554)
(3.456200000000, -0.000053438296)
(3.770400000000, -0.000061174573)
(4.084600000000, -0.000059467472)
(4.398800000000, -0.000052215443)
(4.713000000000, -0.000067396396)
(5.027200000000, -0.000054563621)
(5.341400000000, -0.000059860257)
(5.655600000000, -0.000058522936)
(5.969800000000, -0.000053299221)
(6.284000000000, -0.000067055977)
};
\addlegendentry{\large{$\ell_{\max}=7$, $N_{\leb}=86$}}
\addplot[color=dark_green,  mark=oplus*, line width = 0.5mm, dashdotted,, mark options = {scale= 1.5, solid}] 
coordinates{
(0.000000000000, -0.000056638670)
(0.314200000000, -0.000055413555)
(0.628400000000, -0.000056928516)
(0.942600000000, -0.000056379985)
(1.256800000000, -0.000054663580)
(1.571000000000, -0.000055260420)
(1.885200000000, -0.000054822178)
(2.199400000000, -0.000056248156)
(2.513600000000, -0.000056811287)
(2.827800000000, -0.000055409913)
(3.142000000000, -0.000056616393)
(3.456200000000, -0.000055404332)
(3.770400000000, -0.000057374382)
(4.084600000000, -0.000057363083)
(4.398800000000, -0.000055310289)
(4.713000000000, -0.000056720933)
(5.027200000000, -0.000055392612)
(5.341400000000, -0.000057548465)
(5.655600000000, -0.000057282337)
(5.969800000000, -0.000055714283)
(6.284000000000, -0.000056614652)
};
\addlegendentry{\large{$\ell_{\max}=11$, $N_{\leb}=1202$}}
\end{axis}
\end{tikzpicture}
\caption{Example~\ref{subsec:rot_sym}: Electrostatic solvation energy of the hydrogen-fluoride molecule with respect to the angle of rotating fluorine atom.  
}\label{fig:hf_rot_symmetry}
\end{figure}
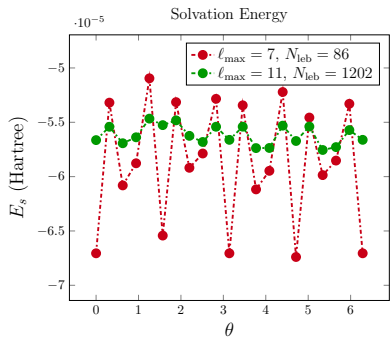

In Fig.~\ref{fig:hf_rot_symmetry}, we plot $E_s$ with respect to $\theta\in [0,2 \pi]$ that is the angle between the atom and notice that variation is systematically controlled with a variance of $6\%$ for $\ell_{\max}=7$ and $2\%$ for $\ell_{\max}=11$. Hence, the variation of the energy under rotational symmetry is systematically controlled as it decreases with increased spherical harmonics.

\subsection{Variation of Debye-H\"uckel Screening Constant}\label{subsec:var_kappa}
In remark~\ref{rem:limiting_case}, we noticed that the ddPCM-SES method can be regarded as a limiting case of ddNPB-SES if $\kappa\rightarrow 0$. For this example, we study the effect of $\kappa$ on the electrostatic solvation energy. We consider the hydrogen fluoride molecule with $r_0=0$\ \AA, and $a=0$\ \AA; and set the discretization parameters as $\ell_{\max}=7$, $N_{\leb}=86$, $N=15$, and $N_{\lgl}=30$.

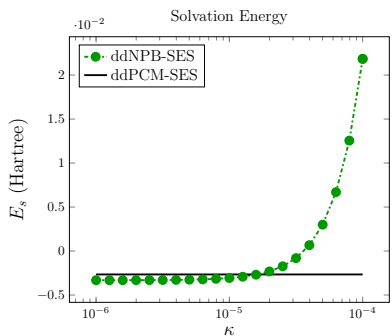
\begin{figure}[!t]
\centering
\begin{tikzpicture}[scale=0.5, spy using outlines={rectangle, magnification=5,connect spies}]
\begin{semilogxaxis}[
    legend pos=north west, xlabel = \Large{$\kappa$}, ylabel= \Large{$E_s$ (Hartree)},  legend cell align ={left}, title = \large{Solvation Energy}]
\addplot[color=dark_green,  mark=oplus*, line width = 0.5mm, dashdotted,,mark options = {scale= 1.5, solid}]
coordinates{
(0.000100000000, 0.021839744316)
(0.000100000000, 0.021839744316)
(0.000079432823, 0.012551995620)
(0.000063095734, 0.006693435130)
(0.000050118723, 0.002997741630)
(0.000039810717, 0.000666322166)
(0.000031622777, -0.000804501460)
(0.000025118864, -0.001732427364)
(0.000019952623, -0.002317857336)
(0.000015848932, -0.002687213119)
(0.000012589254, -0.002920248404)
(0.000010000000, -0.003067276952)
(0.000007943282, -0.003160042805)
(0.000006309573, -0.003218573853)
(0.000005011872, -0.003255500461)
(0.000003981072, -0.003278800400)
(0.000003162278, -0.003293501465)
(0.000002511886, -0.003302777108)
(0.000001995262, -0.003308629592)
(0.000001584893, -0.003312322234)
(0.000001258925, -0.003314652121)
(0.000001000000, -0.003316122174)
};
\addlegendentry{\large{ddNPB-SES}} 
\addplot[color=black,  line width = 0.5mm] 
coordinates{( 0.000001000000,   -0.0026606)( 0.00010000000000,  -0.0026606)};
\addlegendentry{\large{ddPCM-SES}} 
\end{semilogxaxis}
\end{tikzpicture}
\caption{Example~\ref{subsec:var_kappa}:  Electrostatic solvation energy with respect to $\kappa$ for the hydrogen fluoride molecule. The ddPCM-SES energy is obtained with $r_0=0$\ \AA\ in \cite{QSM18}.
}\label{fig:hf_var_kappa}
\end{figure}

We vary $\kappa$ from $10^{-6}$ to $10^{-4}$.  On a continuous level, the ddNPB-SES model converges to the ddPCM-SES model when $\kappa\rightarrow 0$ (see \cite{QSM18}). We notice similar observations numerically; see Fig.~\ref{fig:hf_var_kappa}. As $\kappa$ decreases, the energy becomes constant and converges to a value. We obtain a reference ddPCM-SES result from \cite{QSM18} using $r_0=0$\ \AA.

\section{Summary and Outlook}\label{sec:summary}
This paper proposes a new method for solving the nonlinear Poisson-Boltzmann equation using the domain decomposition paradigm for the solvent-excluded surface.

The original problem defined in $\mathbb{R}^3$ is transformed into two coupled equations described in the bounded solute cavity based on potential theory arguments. An enlarged cavity was defined for each atom encompassing the Stern layer and a nonlinear regime. Then, the Schwarz domain decomposition method was used to solve these two problems by decomposing them into balls. We develop two single-domain solvers for solving the HSP and GSP equations in a unit ball. The GSP solver was a nonlinear solver that used spherical harmonics for angular direction and Legendre polynomials for radial direction. An SES-based continuous dielectric permittivity function and an ion-exclusion function were proposed and were encompassed in the GSP solver. The novelty of the method is that the nonlinearity is incorporated only in the proximity of the molecule, whereas in the bulk solvent region, the linear model was used. A series of numerical results have been presented to show the performance of the ddNPB-SES method, which shows the importance of using the nonlinear regime near the surface.

In the future, we would like to implement this method in our open-source software \texttt{ddX}, \cite{ddX}, to simulate bigger molecules. We would also like to study the nonlinear solver better using Newton's methods and incorporating acceleration techniques such as dynamic damping and Anderson acceleration.

\bibliographystyle{plain}
\bibliography{ddNLPB}
\end{document}